\documentclass[mathpazo]{cicp}

\usepackage{color}
\newcommand{\Q}{\mathbf{Q}}
\newcommand{\F}{\mathbf{F}}
\renewcommand{\S}{\mathbf{S}}
\newcommand{\K}{\mathbf{K}}
\newcommand{\M}{\mathbf{M}}
\renewcommand{\P}{\mathbf{P}}
\renewcommand{\v}{\mathbf{v}}
\newcommand{\w}{\mathbf{w}}
\newcommand{\x}{\mathbf{x}}

\newcommand{\q}{\mathbf{q}}
\newcommand{\f}{\mathbf{f}}
\newcommand{\g}{\mathbf{g}}

\begin{document}
\title{Arbitrary-Lagrangian-Eulerian One-Step WENO \mbox{Finite} Volume Schemes on Unstructured Triangular Meshes}

\author[Dumbser et.~al.]{Walter Boscheri\affil{1},
                         Michael Dumbser\affil{1}\comma\corrauth} 
\address{\affilnum{1}\ Laboratory of Applied Mathematics, Department of Civil, Environmental and \\ Mechanical Engineering, 
          					   University of Trento,
          						 I-38123 Trento, Italy
          				     }
\emails{{\tt walter.boscheri@unitn.it} (W.~Boscheri), \\ {\tt michael.dumbser@unitn.it} (M.~Dumbser)}

\begin{abstract}
In this article we present a new class of high order accurate Arbitrary--Eulerian--Lagrangian (ALE) one--step WENO finite volume schemes for solving nonlinear hyperbolic systems of conservation laws 
on moving two dimensional unstructured triangular meshes. A WENO reconstruction algorithm is used to achieve high order accuracy in space and a high order one--step time discretization is achieved by using 
the local space--time Galerkin predictor proposed in \cite{Dumbser20088209}. For that purpose, a new element--local weak formulation of the governing PDE is adopted on moving space--time elements. 
The space-–time basis and test functions are obtained considering Lagrange interpolation polynomials passing through a predefined set of nodes. Moreover, a polynomial mapping 
defined by the same local space--time basis functions as the weak solution of the PDE is used to map the moving physical space--time element onto a space--time reference element. To maintain 
algorithmic simplicity, the final ALE one--step finite volume scheme uses moving triangular meshes with \textit{straight} edges. This is possible in the ALE framework, which allows a local mesh 
velocity that is different from the local fluid velocity. We present numerical convergence rates for the schemes presented in this paper up to sixth order of accuracy in space and time and show 
some classical numerical test problems for the two--dimensional Euler equations of compressible gas dynamics.   
\end{abstract}


\maketitle

\section{Introduction}
\label{sec.intro} 
In this paper we present a new family of high order accurate Lagrangian--type one--step finite volume schemes for solving nonlinear hyperbolic balance laws, with non stiff algebraic source term. The main advantage of working in 
a Lagrangian framework is that material interfaces can be identified and located precisely, since the computational mesh is moving with the local fluid velocity, hence obtaining a more accurate resolution of material  interfaces. For this reason a lot of research has been carried out in the last decades in order to develop Lagrangian methods, whose algorithms can start either directly from the conservative quantities such as mass, momentum 
and total energy \cite{Maire2007,Smith1999}, or from the nonconservative form of the governing equations, as proposed in \cite{Benson1992,Caramana1998,Neumann1950}. Furthermore we can split the existing Lagrangian schemes into 
two main classes, depending on the location of the physical variables on the mesh: in one case the velocity is defined at the cell interfaces and the other variables at the cell barycenter, hence adopting a \textit{staggered mesh}  
approach, while in the other case all variables are defined at the cell barycenter, therefore using a \textit{cell--centered} approach.  

The equations of Lagrangian gas dynamics have been considered in \cite{munz94}, where several different Godunov-type finite volume schemes have been presented and where a new Roe linearization has been introduced in order to 
define proper estimates of the maximum signal speeds in HLL--type Riemann solvers in Lagrangian coordinates. A cell-centered Godunov scheme has been proposed by  Carr\'e et al. \cite{Després2009} for Lagrangian gas dynamics 
on general multi-dimensional unstructured meshes. In this case the finite volume scheme is node based and compatible with the mesh displacement. In \cite{DepresMazeran2003} Despr\'es and Mazeran introduce a new formulation of 
the multidimensional Euler equations in Lagrangian coordinates as a system of conservation laws associated with constraints. Furthermore they propose a way to evolve in a coupled manner both the physical and the geometrical 
part of the system \cite{Despres2005}, writing the two--dimensional equations of gas dynamics in Lagrangian coordinates together with the evolution of the geometry as a weakly hyperbolic system of conservation laws. This allows 
the authors to design a finite volume scheme for the discretization of Lagrangian gas dynamics on moving meshes, based on the symmetrization of the formulation of the physical part. In a recent work Despr\'es et al.  \cite{Depres2012} propose a new method designed for cell-centered Lagrangian schemes, which is translation invariant and suitable for curved meshes. General polygonal grids are also considered by Maire et al. \cite{Maire2009,Maire2010,Maire2011}, who develop a general formalism to derive first and second order cell-centered Lagrangian schemes in multiple space dimensions. By the use of a node-centered solver \cite{Maire2009}, 
the authors obtain the time derivatives of the fluxes. The solver may be considered as a multi-dimensional extension of the Generalized Riemann problem methodology introduced by Ben-Artzi and Falcovitz \cite{Artzi}, Le Floch 
et al. \cite{Raviart.GRP.1,Raviart.GRP.2} and Titarev and Toro \cite{Toro:2006a,toro3,toro10}. So far, all the above--mentioned schemes are at most second order accurate in space and time.  

In order to achieve higher accuracy, Cheng and Shu were the first who introduced a high order essentially non-oscillatory (ENO) reconstruction in Lagrangian schemes \cite{chengshu1,chengshu2}. They developed a class of cell centered Lagrangian finite volume schemes for solving the Euler equations, using both Runge-Kutta and Lax-Wendroff-type time stepping to achieve also higher order in time. Furthermore a formalism for symmetry preserving 
Lagrangian schemes has been proposed by Cheng and Shu, see \cite{chengshu3,chengshu4}. The higher order schemes presented in \cite{chengshu1,chengshu2} were the first better than second order non--oscillatory Lagrangian--type 
finite volume schemes ever proposed. Higher order finite element methods on unstructured meshes have been investigated in \cite{scovazzi1,scovazzi2}. In a very recent paper Dumbser et al. \cite{Dumbser2012} propose a new 
class of high order accurate Lagrangian--type one-–step WENO finite volume schemes for the solution of stiff hyperbolic balance laws. They consider the one--dimensional case and develop a Lagrangian algorithm up to eighth 
order of accuracy, based on high order WENO reconstruction in space and the local space--time discontinuous Galerkin predictor method proposed in \cite{DumbserEnauxToro} to obtain a high order one--step scheme in time. 

Other schemes can be also included and mentioned as Lagrangian algorithms, e.g. meshless particle schemes that adopt a fully Lagrangian approach, such as the smooth particle hydrodynamics (SPH) method \cite{Monaghan1994,SPHLagrange,SPHWeirFlow,SPH3D,Dambreak3D}, which can be used to simulate fluid motion in complex deforming domains. Also within the SPH approach, which is a fully Lagrangian method, the mesh moves with 
the local fluid velocity, whereas in Arbitrary Lagrangian Eulerian (ALE) schemes, see e.g. \cite{Hirt1974,Peery2000,Smith1999,Feistauer1,Feistauer2,Feistauer3,Feistauer4}, the mesh moves with an arbitrary mesh velocity that does 
not necessarily coincide with the real fluid velocity. Furthermore one can find Semi-Lagrangian schemes, which are mainly used for solving transport equations \cite{ALE2000Belgium,ALE1996FV}. Here, the numerical solution at 
the new time level is computed from the known solution at the present time by following backward in time the Lagrangian trajectories of the fluid to the end-point of the trajectory. Since the end-point does not coincide with 
a grid point, an interpolation formula is required in order to evaluate the unknown solution, see e.g.  \cite{Casulli1990,CasulliCheng1992,LentineEtAl2011,HuangQiul2011,QuiShu2011,CIR,BoscheriDumbser}. For the sake of clarity 
we specify that in Semi-Lagrangian algorithms the mesh is fixed as in a classical Eulerian approach. An alternative to Lagrangian methods for the accurate resolution of material interfaces has been developed in the Eulerian framework on fixed meshes under the form of the ghost--fluid method \cite{FedkiwEtAl1,FedkiwEtAl2,FerrariLevelSet}, together with a level--set approach \cite{levelset1,levelset2}, where the level--set function represents the 
signed distance from the material interface and its zeros locate the interface position.  
 
In this paper we introduce a new better than second order accurate \textit{two--dimensional} Lagrangian--type one--step WENO finite volume scheme on \textit{unstructured triangular meshes}: 
high order of accuracy in space is obtained using a WENO reconstruction \cite{balsara,JiangShu1996,Dumbser2007693,Dumbser2007204,MixedWENO2D,MixedWENO3D,HuShuVortex1999,ZhangShu3D,AboiyarIske}, although other 
higher order spatial reconstruction schemes could be adopted as well, see \cite{abgrall_eno,MOOD}. High order accuracy in time is achieved with a local space--time Galerkin predictor, as introduced in  \cite{DumbserEnauxToro,DumbserZanotti,HidalgoDumbser,Dumbser2012}. The method proposed in this article is presented as an Arbitrary--Lagrangian--Eulerian scheme in order to allow arbitrary grid motion. This allows us 
to use curved space--time elements in the local predictor stage but triangles delimited by straight line segments in the resulting one--step finite volume scheme (corrector stage). This choice has been made to maintain 
algorithmic simplicity. 
 
The outline of this article is as follows: in Section \ref{sec.method} we present the details of the proposed numerical scheme, while in Section \ref{sec.test} we show numerical convergence rates up to sixth order of 
accuracy in space and time for a smooth problem as well as numerical results for several different test cases governed by the compressible Euler equations. The paper closes with some concluding remarks and an 
outlook to possible future extensions of the method in Section \ref{sec.conclusions}. 

\section{Numerical Method}
\label{sec.method}

In this article we consider general nonlinear systems of hyperbolic balance laws of the form
\begin{equation}
\label{PDE}
  \frac{\partial \Q}{\partial t} + \nabla \cdot \F(\Q) = \S(\Q), \qquad (x,y) \in \Omega(t) \subset \mathbb{R}^2, \quad t \in \mathbb{R}_0^+, \quad \Q \in \Omega_{\Q} \subset \mathbb{R}^\nu,     
\end{equation} 
where $\Q=(q_1,q_2,...,q_\nu)$ is the vector of conserved variables defined in the space of the admissible states $\Omega_{\Q} \subset \mathbb{R}^\nu$, $\F(\Q)=\left( \f(\Q),\g(\Q) \right)$ is the nonlinear flux tensor and $\S(\Q)$ 
represents a nonlinear but non-stiff algebraic source term. The spatial position vector is denoted by $\mathbf{x}=(x,y)$, $t$ is the time and $\Omega(t)$ is the time--dependent computational domain. 

The two-dimensional time--dependent computational domain $\Omega(t)$ is discretized at the current time $t^n$ by a set of triangular elements $T^n_i$. The union of all elements is called the 
\textit{current triangulation} $\mathcal{T}^n_{\Omega}$ of the domain $\Omega(t^n)=\Omega^n$ and can be expressed as  
\begin{equation}
\mathcal{T}^n_{\Omega} = \bigcup \limits_{i=1}^{N_E}{T^n_i}, 
\label{trian}
\end{equation}
where $N_E$ is the total number of elements contained in the domain.  

In a Lagrangian framework we are dealing with a \textit{moving} mesh, hence the elements deform while the solution is evolving in time. It is therefore convenient 
to adopt a \textit{local} reference coordinate system $\xi-\eta$, where the reference element $T_e$ is defined. The spatial mapping of the triangular 
elements $T^n_i$ at the current time $t^n$ from reference coordinates $\xi-\eta$ to physical coordinates $x-y$ is given by the relation  
\begin{eqnarray} 
 x &=& X^n_{1,i} + \left( X^n_{2,i} - X^n_{1,i} \right) \xi + \left( X^n_{3,i} - X^n_{1,i} \right) \eta, \nonumber \\ 
 y &=& Y^n_{1,i} + \left( Y^n_{2,i} - Y^n_{1,i} \right) \xi + \left( Y^n_{3,i} - Y^n_{1,i} \right) \eta,    
 \label{xietaTransf} 
\end{eqnarray} 
where $\mathbf{X}^n_{k,i} = (X^n_{k,i},Y^n_{k,i})$ denotes the vector of spatial coordinates of the $k$-th vertex of the triangle $T^n_i$ in physical coordinates at the current time $t^n$. 
The vector $\boldsymbol{\xi} = (\xi, \eta)$ is the vector of spatial coordinates in the reference system, while $\mathbf{x}=(x,y)$ is the spatial coordinate vector in the physical 
system. The spatial reference element $T_e$ is the unit triangle composed of the nodes $\boldsymbol{\xi}_{e,1}=(\xi_{e,1},\eta_{e,1})=(0,0)$, $\boldsymbol{\xi}_{e,2}=(\xi_{e,2},\eta_{e,2})=(1,0)$ 
and $\boldsymbol{\xi}_{e,3}=(\xi_{e,2},\eta_{e,2})=(0,1)$.  

As usual for finite volume schemes, data are represented by spatial cell averages, which are defined at time $t^n$ as  
\begin{equation}
  \Q_i^n = \frac{1}{|T_i^n|} \int_{T^n_i} \Q(x,y,t^n) dV,     
 \label{eqn.cellaverage}
\end{equation}  
where $|T_i^n|$ denotes the volume of element $T_i^n$ at the current time $t^n$. To achieve higher order in space, piecewise high order polynomials $\mathbf{w}_h(x,y,t^n)$ must be 
reconstructed from the given cell averages \eqref{eqn.cellaverage} using a polynomial WENO reconstruction procedure illustrated briefly in the next section.  

\subsection{Polynomial WENO Reconstruction on Unstructured Meshes}
\label{sec.WENOrec}

In this paper we use the WENO reconstruction algorithm in the \textit{polynomial} formulation presented in \cite{DumbserEnauxToro,Dumbser2007204,Dumbser2007693}, instead of adopting the original \textit{pointwise} WENO scheme \cite{JiangShu1996,HuShuVortex1999,ZhangShu3D}. Here we will give a very brief summary of the algorithm, since it is completely described in all detail in the above-mentioned references.  

The reconstruction is performed by using the \textit{reference system} $(\xi,\eta)$ according to the mapping \eqref{xietaTransf}, as also explained in detail in \cite{Dumbser2007693}. 
The reconstruction polynomial of degree $M$ is obtained by considering a reconstruction stencil $S_i^s$ 
\begin{equation}
\mathcal{S}_i^s = \bigcup \limits_{j=1}^{n_e} T^n_{m(j)}. 
\label{stencil}
\end{equation}
The stencil contains a total number of elements $n_e$ that is greater than the smallest number $\mathcal{M} = (M+1)(M+2)/2$ needed to reach the formal order of accuracy $M+1$, as shown in  
\cite{StencilRec1990,Olliver2002,KaeserIske2005}. Typically we take $n_e = 2 \mathcal{M}$ in two space dimensions. Here $1\leq j \leq n_e$ is a local index counting the elements belonging 
to the stencil, while $m(j)$ maps the local index to the global element numbers used in the triangulation \eqref{trian}.  
$s$ denotes the number and the position of the stencil with respect to the central element $T^n_i$. According to \cite{Dumbser2007693,Dumbser2007204} we always will use one central 
reconstruction stencil given by $s=0$, three  primary sector stencils and three reverse sector stencil, as suggested by K\"aser and Iske \cite{KaeserIske2005}, so that we globally deal with 
seven stencils per element.    

We use some spatial basis functions $\psi_l(\xi,\eta)$ in order to write the reconstruction polynomial for each candidate stencil $s$ for triangle $T_i^n$: 
\begin{equation}
\label{eqn.recpolydef} 
\w^s_h(x,y,t^n) = \sum \limits_{l=1}^\mathcal{M} \psi_l(\xi,\eta) \hat \w^{n,s}_{l,i} := \psi_l(\xi,\eta) \hat \w^{n,s}_{l,i},   
\end{equation}
with the mapping to the reference coordinate system given by the transformation \eqref{xietaTransf}. In the rest of the paper we will use classical tensor index notation with the Einstein summation
convention, which implies summation over two equal indices. The number of the unknown polynomial coefficients (degrees of freedom) is $\mathcal{M}$, already defined previously. As basis functions 
$\psi_l(\xi,\eta)$ we use the orthogonal Dubiner--type basis on the reference triangle $T_e$, given in detail in \cite{Dubiner,orth-basis,CBS-book}. 

Integral conservation is required for the reconstruction on each element $T^n_j \in \mathcal{S}_i^s$, hence  
\begin{equation}
\label{intConsRec}
\frac{1}{|T^n_j|} \int \limits_{T^n_j} \psi_l(\xi,\eta) \hat \w^{n,s}_{l,i} dV = \Q^n_j, \qquad \forall T^n_j \in \mathcal{S}_i^s     
\end{equation}
where $|T^n_j|$ represents the volume of element $T^n_j$ at time $t^n$. Since the number of stencil elements is \textit{larger} than the one of the unknown polynomial coefficients ($n_e > \mathcal{M}$), the 
above system given by eqn. \eqref{intConsRec} is an overdetermined linear algebraic system that is readily solved for the unknown coefficients $\hat \w^{n,s}_{l,i}$ using a constrained least--squares technique, 
see \cite{Dumbser2007693}.  
The multi--dimensional integrals are evaluated using Gaussian quadrature formulae of suitable order, see the book of Stroud for details \cite{stroud}. 
Since the shape of the triangles changes in time, the small linear systems given by \eqref{intConsRec} must be solved at the beginning of each time step, while the choice of the stencils $\mathcal{S}_i^s$ remains 
\textit{fixed} for all times. It is therefore very useful to devise a one--step time--discretization, such as the one detailed in the next section, which requires only one reconstruction per time step, in contrast 
to higher order Runge--Kutta methods, which would require the evaluation of the above reconstruction operator in each substage of the Runge--Kutta method.  

In order to obtain essentially non--oscillatory properties the final WENO reconstruction polynomial is computed from the reconstruction polynomials obtained on each individual stencil $\mathcal{S}_i^s$ 
in the usual way. For this purpose we adopt the classical definitions of the oscillation indicators $\sigma_s$  
reported in \cite{JiangShu1996} and the oscillation indicator matrix $\Sigma_{lm}$ proposed in \cite{Dumbser2007204,Dumbser2007693}, which read  
\begin{equation}
\sigma_s = \Sigma_{lm} \hat w^{n,s}_{l,i} \hat w^{n,s}_{m,i}, \qquad 
\Sigma_{lm} = \sum \limits_{ \alpha + \beta \leq M}  \, \, \int \limits_{T_e} \frac{\partial^{\alpha+\beta} \psi_l(\xi,\eta)}{\partial \xi^\alpha \partial \eta^\beta} \cdot 
                                                                         \frac{\partial^{\alpha+\beta} \psi_m(\xi,\eta)}{\partial \xi^\alpha \partial \eta^\beta} d\xi d\eta.   
\end{equation} 
The nonlinear weights $\omega_s$ are defined by
\begin{equation}
\tilde{\omega}_s = \frac{\lambda_s}{\left(\sigma_s + \epsilon \right)^r}, \qquad 
\omega_s = \frac{\tilde{\omega}_s}{\sum_q \tilde{\omega}_q},  
\end{equation} 
where we use $\epsilon=10^{-14}$, $r=8$, $\lambda_s=1$ for the one--sided stencils and $\lambda_0=10^5$ for the central stencil, according to \cite{DumbserEnauxToro,Dumbser2007204}. 
The final nonlinear WENO reconstruction polynomial and its coefficients are then given by 
\begin{equation}
\label{eqn.weno} 
 \w_h(x,y,t^n) = \sum \limits_{l=1}^{\mathcal{M}} \psi_l(\xi,\eta) \hat \w^{n}_{l,i}, \qquad \textnormal{ with } \qquad  
 \hat \w^{n}_{l,i} = \sum_s \omega_s \hat \w^{n,s}_{l,i}.   
\end{equation}

\subsection{Local Space--Time Continuous--Galerkin Predictor on Moving Meshes}
\label{sec.localCG}

In order to obtain high order accuracy in time, the reconstructed polynomials $\w_h$ obtained at the current time $t^n$ are now evolved \textit{locally} within each element $T_i(t)$ during one 
time step $t^n \leq t \leq t^{n+1}$. The result of this local evolution step is an \textit{element--local} piecewise polynomial space--time predictor $\q_h(x,y,t)$. For this purpose we use an 
element--local weak formulation of the PDE in both space and time, based on a Lagrangian version of the local space--time continuous Galerkin method introduced for the Eulerian framework in 
\cite{Dumbser20088209}. 

Let $\mathbf{\tilde{x}}=(x,y,t)$ denote the physical coordinate vector and $\boldsymbol{\tilde{\xi}}=(\xi,\eta,\tau)$ the reference coordinate vector, in both of which also the time is included, 
while $\mathbf{x}=(x,y)$ and $\boldsymbol{\xi}=(\xi,\eta)$ are the pure spatial coordinate vectors in physical and reference coordinates, respectively. 
Let furthermore $\theta_l=\theta_l(\boldsymbol{\tilde{\xi}})=\theta_l(\xi,\eta,\tau)$ be a space--time basis function defined by the Lagrange interpolation polynomials passing through the space--time 
nodes $\boldsymbol{\tilde{\xi}}_m=(\xi_m,\eta_m,\tau_m)$. The nodes are specified according to \cite{Dumbser20088209}. Being the Lagrange interpolation polynomials, the resulting \textit{nodal} basis 
functions $\theta_l$ therefore satisfy the interpolation property 
\begin{equation}
 \theta_l(\boldsymbol{\tilde{\xi}}_m) = \delta_{lm}, 
\end{equation} 
where $\delta_{lm}$ is the usual Kronecker symbol. 
According to \cite{Dumbser20088209} the local space–-time solution in element $T_i(t)$, denoted by $\q_h$, as well as the fluxes $(\f_{h}, \g_h)$ and the source term $\S_h$ are approximated as 
\begin{eqnarray}
\q_h=\q_h(\xi,\eta,\tau) = \theta_{l}(\xi,\eta,\tau) \widehat{\q}_{l,i}, \qquad & \S_h=\S_h(\xi,\eta,\tau) = \theta_{l}(\xi,\eta,\tau) \widehat{\S}_{l,i},  \nonumber\\
\f_h=\f_h(\xi,\eta,\tau) = \theta_{l}(\xi,\eta,\tau) \widehat{\f}_{l,i}, \qquad & \g_h=\g_h(\xi,\eta,\tau) = \theta_{l}(\xi,\eta,\tau) \widehat{\g}_{l,i}.
\label{thetaSol}
\end{eqnarray}
The mapping between the physical space--time coordinate vector $\mathbf{\tilde{x}}$ and the reference space--time coordinate vector $\boldsymbol{\tilde{\xi}}$ is achieved using an \textit{isoparametric} 
approach, i.e. the mapping is represented by the \textit{same} basis functions $\theta_l$ as the approximate numerical solution itself, given in Eqn. \eqref{thetaSol} above. Hence, one has 
\begin{equation}
 x(\xi,\eta,\tau) = \theta_l(\xi,\eta,\tau) \widehat{x}_{l,i}, \qquad y(\xi,\eta,\tau) = \theta_l(\xi,\eta,\tau) \widehat{y}_{l,i}, \qquad t(\xi,\eta,\tau) = \theta_l(\xi,\eta,\tau) \widehat{t}_l, 
 \label{eqn.isoparametric} 
\end{equation} 
where the degrees of freedom $\widehat{\mathbf{x}}_{l,i} = (\widehat{x}_{l,i},\widehat{y}_{l,i})$ denote the (in parts unknown) vector of physical coordinates in space of the moving space--time control volume 
and the $\widehat{t}_l$ denote the \textit{known} degrees of freedom of the physical time at each space--time node $\tilde{\x}_{l,i} = (\widehat{x}_{l,i}, \widehat{y}_{l,i}, \widehat{t}_l)$.  
The mapping in time simply reads  
\begin{equation}
t = t_n + \tau \, \Delta t, \qquad  \tau = \frac{t - t^n}{\Delta t}, \qquad \Rightarrow \qquad \widehat{t}_l = t_n + \tau_l \, \Delta t, 
\label{timeTransf}
\end{equation} 
where $t^n$ is the current time and $\Delta t$ is the time step. Hence, $t_\xi = t_\eta = 0$ and $t_\tau = \Delta t$. 
    
The spatial mapping given in (\ref{eqn.isoparametric}) and the temporal mapping \eqref{timeTransf} allow us to transform the physical space--time element to the unit reference space--time element 
$T_e \times [0,1]$. The Jacobian of the spatial and temporal transformation is given by 
\begin{equation}
J_{st} = \frac{\partial \mathbf{\tilde{x}}}{\partial \boldsymbol{\tilde{\xi}}} = \left( \begin{array}{ccc} x_{\xi} & x_{\eta} & x_{\tau} \\ y_{\xi} & y_{\eta} & y_{\tau} \\ 0 & 0 & \Delta_t \\ \end{array} \right), 
\label{Jac}
\end{equation}
and its inverse reads
\begin{equation}
J_{st}^{-1} = \frac{\partial \boldsymbol{\tilde{\xi}}}{\partial \mathbf{\tilde{x}}} = \left( \begin{array}{ccc} \xi_{x} & \xi_{y} & \xi_{t} \\ \eta_{x} & \eta_{y} & \eta_{t} \\ 0 & 0 & \frac{1}{\Delta t} \\ \end{array} \right).
\label{iJac}
\end{equation}

The local reference system and the inverse of the associated Jacobian matrix (\ref{iJac}) are used to rewrite the governing PDE (\ref{PDE}) as 
\begin{equation}
\frac{\partial \Q}{\partial \tau}\tau_t + \frac{\partial \Q}{\partial \xi}\xi_t + \frac{\partial \Q}{\partial \eta}\eta_t + \frac{\partial \f}{\partial \tau}\tau_x + \frac{\partial \f}{\partial \xi}\xi_x + \frac{\partial \f}{\partial \eta}\eta_x + \frac{\partial \g}{\partial \tau}\tau_y + \frac{\partial \g}{\partial \xi}\xi_y + \frac{\partial \g}{\partial \eta}\eta_y = \mathbf{S}(\Q),  
\label{PDEweak}
\end{equation}
which simplifies to 
\begin{equation}
\frac{\partial \Q}{\partial \tau} + \Delta t \left[ \frac{\partial \Q}{\partial \xi}\xi_t + \frac{\partial \Q}{\partial \eta}\eta_t + \frac{\partial \f}{\partial \xi}\xi_x + \frac{\partial \f}{\partial \eta}\eta_x + \frac{\partial \g}{\partial \xi}\xi_y + \frac{\partial \g}{\partial \eta}\eta_y  \right] = \Delta t \mathbf{S}(\Q),
\label{PDECG}
\end{equation}
since $\tau_x = \tau_y = 0$ and $\tau_t = \frac{1}{\Delta t}$, according to the definition (\ref{timeTransf}).

To make the notation easier we now introduce the following two integral operators
\begin{equation}
\left[f,g\right]^{\tau} = \int \limits_{T_e} f(\xi,\eta,\tau) g(\xi,\eta,\tau) d\xi d\eta, \quad \left\langle f,g \right\rangle = \int \limits_{0}^{1} \int \limits_{T_e} f(\xi,\eta,\tau)g(\xi,\eta,\tau) d\xi d\eta d\tau,  
\label{intOperators}
\end{equation}
which denote the scalar products of two functions $f$ and $g$ over the spatial reference element $T_e$ at time $\tau$ and over the space-time reference element $T_e\times \left[0,1\right]$, respectively.

Inserting the definitions of (\ref{thetaSol}) into the weak formulation of the PDE (\ref{PDECG}), we multiply Eqn. \eqref{PDECG} with the same space--time basis functions $\theta_k(\xi,\eta,\tau)$ and then integrate it over the space--time reference element $T_e \times [0,1]$, hence obtaining:
\begin{eqnarray}
\left\langle \theta_k,\frac{\partial \theta_l}{\partial \tau} \right\rangle \widehat{\q}_{l,i} 
&+& \Delta t \left( \left\langle \theta_k,\frac{\partial \theta_l}{\partial \xi}\xi_t \right\rangle + \left\langle \theta_k,\frac{\partial \theta_l}{\partial \eta}\eta_t \right\rangle \right) \widehat{\q}_{l,i} \nonumber\\
&+& \Delta t \left( \left\langle \theta_k,\frac{\partial \theta_l}{\partial \xi}\xi_x \right\rangle + \left\langle \theta_k,\frac{\partial \theta_l}{\partial \eta}\eta_x \right\rangle \right) \widehat{\f}_{l,i} \nonumber\\
&+& \Delta t \left( \left\langle \theta_k,\frac{\partial \theta_l}{\partial \xi}\xi_y \right\rangle + \left\langle \theta_k,\frac{\partial \theta_l}{\partial \eta}\eta_y \right\rangle \right) \widehat{\g}_{l,i} 
= \left\langle \theta_k,\theta_l \right\rangle \widehat{\S}_{l,i}. \nonumber\\ 
\label{LagrSTPDECG}
\end{eqnarray}

The above expression can be written in a more compact matrix form, which reads
\begin{equation}
\K_{\tau}\widehat{\q}_{l,i} + \Delta t \left( \K_t \widehat{\q}_{l,i} + \K_x \widehat{\f}_{l,i} + \K_y \widehat{\g}_{l,i} \right),  
= \Delta t \M \widehat{\S}_{l,i}, 
\label{LagrSTPDECGmatrix}
\end{equation}
with the following matrix definitions:
\begin{equation}
\K_{\tau} = \left\langle \theta_k,\frac{\partial \theta_l}{\partial \tau} \right\rangle, \qquad \M = \left\langle \theta_k,\theta_l \right\rangle, 
\label{Ktau}
\end{equation}
and 
\begin{eqnarray}
\K_t  &=& \left\langle \theta_k,\frac{\partial \theta_l}{\partial \xi}\xi_t \right\rangle + \left\langle \theta_k,\frac{\partial \theta_l}{\partial \eta}\eta_t \right\rangle,  \nonumber \\
\K_x  &=& \left\langle \theta_k,\frac{\partial \theta_l}{\partial \xi}\xi_x \right\rangle + \left\langle \theta_k,\frac{\partial \theta_l}{\partial \eta}\eta_x \right\rangle,  \nonumber \\ 
\K_y  &=& \left\langle \theta_k,\frac{\partial \theta_l}{\partial \xi}\xi_y \right\rangle + \left\langle \theta_k,\frac{\partial \theta_l}{\partial \eta}\eta_y \right\rangle. 
\label{matricesPDECG}
\end{eqnarray}
Note that the Lagrangian nature of the scheme, i.e. the moving space--time control volume, leads to the term $\K_t \widehat{\q}_h$, which is not present in the Eulerian case introduced 
in \cite{Dumbser20088209}. 
In \eqref{LagrSTPDECGmatrix}, only the matrix $\K_\tau$ and the space--time mass matrix $\mathbf{M}$ given in (\ref{Ktau}) can be pre--computed and stored. The other integrals continuously change 
in time, since they depend on the inverse of the Jacobian matrix, which varies according to the Lagrangian mesh motion. Therefore it may be more convenient to assemble a unified term 
$\P$ as  
\begin{equation}
 \P := \S(\Q) - \left( \Q_{\xi}\xi_t + \Q_{\eta}\eta_t + \f_{\xi}\xi_x + \f_{\eta}\eta_x + \g_{\xi}\xi_y + \g_{\eta}\eta_y \right), 
 \label{eqn.pdef}
\end{equation} 
hence 
\begin{equation}
\Q_\tau =  \Delta t \P.
\label{PCG}
\end{equation} 
$\P$ is approximated again inside element $T_i(t)$ by adopting a nodal approach as done in (\ref{thetaSol}), i.e. 
\begin{equation}
\mathbf{P}_h=\mathbf{P}_h(\xi,\eta,\tau) = \theta_{l}(\xi,\eta,\tau) \widehat{\P}_{l,i},
\label{Pterm}
\end{equation}
with $\widehat{\P}_{l,i}= \P(\mathbf{\tilde{x}}_{l,i})$. 

Let us denote with $\widehat{\q}_{l,i}^{0}$ the degrees of freedom that are known from the initial condition $\w_h$ so that $\q_h(x,y,t^n)=\w_h(x,y,t^n)$, whereas the unknown degrees of freedom for 
$\tau>0$ are denoted by $\widehat{\q}_{l,i}^{1}$, so that the total vector of degrees of freedom is written as $\widehat{\q}_{l,i}=(\widehat{\q}_{l,i}^{0},\widehat{\q}_{l,i}^{1})$. 
One can move onto the right-hand side of \eqref{LagrSTPDECGmatrix} or (\ref{PCG}) all the known degrees of freedom $\widehat{\q}_{l,i}^{0}$ given by the 
initial condition $\w_h$ by setting the corresponding degrees of freedom to the known values, like a standard Dirichlet boundary condition in the continuous finite element framework, see 
\cite{Dumbser20088209} for more details. The final expression for the iterative scheme solving the nonlinear algebraic equation system of the two-dimensional Lagrangian continuous Galerkin predictor method 
simply reads  
\begin{equation}
  \K_{\tau} \widehat{\q}_{l,i}^{r+1} = \Delta t \M \widehat{\P}_{l,i}^r,  
\label{CGfinal}
\end{equation}
where the superscript $r$ denotes the iteration number here. For an efficient initial guess ($r=0$) based on a second order MUSCL--type scheme, see \cite{HidalgoDumbser}. 

Since the mesh is moving we also have to consider the evolution of the vertex coordinates of the local space--time element, whose motion is described by the ODE system 
\begin{equation}
\frac{d \mathbf{x}}{dt} = \mathbf{V}(x,y,t),
\label{ODEmesh}
\end{equation}
where $\mathbf{V}=\mathbf{V}(x,y,t)=(U,V)$ is the local mesh velocity. Since we are developing an \textit{arbitrary Lagrangian-Eulerian} scheme (ALE), we want the mesh velocity to be independent from the fluid velocity. 
In this framework we can deal in the same way with both Eulerian and Lagrangian schemes: if $\mathbf{V}=0$ the scheme reduces indeed to a pure Eulerian approach, while if $\mathbf{V}$ coincides with the local fluid 
velocity $\mathbf{v}$ we obtain a Lagrangian--type method. The velocity inside element $T_i(t)$ can be also expressed as 
\begin{equation}
\mathbf{V}_h=\mathbf{V}_h(\xi,\eta,\tau) = \theta_{l}(\xi,\eta,\tau) \widehat{\mathbf{V}}_{l,i}, 
\label{Vdof}
\end{equation}
with $ \widehat{\mathbf{V}}_{l,i} = \mathbf{V}(\mathbf{\tilde{x}}_{l,i})$.

As introduced in \cite{Dumbser2012}, we can solve the ODE (\ref{ODEmesh}) for the unknown coordinate vector $\widehat{\mathbf{x}}_l$ using again the local space--time CG method: 
\begin{equation}
\left\langle \theta_k,\frac{\partial \theta_l}{\partial \tau} \right\rangle \widehat{\mathbf{x}}_{l,i} = \Delta t \left\langle \theta_k,\theta_l \right\rangle \widehat{\mathbf{V}}_{l,i},
\label{VCG}
\end{equation}
hence resulting in the following iteration scheme for the element--local space--time predictor for the nodal coordinates:  
\begin{equation}
\K_{\tau} \widehat{\mathbf{x}}^{r+1}_{l,i} = \Delta t \M \widehat{\mathbf{V}}^r_{l,i}.
\label{newVertPos}
\end{equation}
The nodal degrees of freedom $\widehat{\x}_l$ at relative time $\tau=0$, i.e. the initial condition of the ODE system, are given by the spatial mapping \eqref{xietaTransf}, since the physical triangle $T_i^n$ at time 
$t^n$ is known. 

Eqn. (\ref{VCG}) is iterated \textit{together} with the weak formulation for the solution (\ref{CGfinal}). The algorithm for the two--dimensional local space--time CG predictor for moving meshes can be summarized by the following steps:
\begin{itemize}
	\item compute the local mesh velocity (\ref{Vdof}), usually by choosing the fluid velocity, hence $\widehat{\mathbf{V}}_{l,i} = \mathbf{v}(\tilde{\x}_{l,i})$;
	\item with (\ref{newVertPos}) update the geometry \textit{locally} within the predictor stage, i.e. the element--local space--time coordinates; 
	\item compute the Jacobian matrix and its inverse by using (\ref{Jac})-(\ref{iJac});
	\item compute the term $\P$ according to (\ref{PCG});
	\item evolve the solution \textit{locally} with (\ref{CGfinal}).
\end{itemize}

The iterative procedure described above stops when the residuals of (\ref{CGfinal}) and \eqref{newVertPos} are less than a prescribed tolerance $tol$ (typically $tol\approx 10^{-12}$). 

Once we have carried out the above procedure for all the elements of the computational domain, we end up with an \textit{element--local predictor} for the numerical solution $\q_h$, for the fluxes $\mathbf{F}_h=(\f_h,\g_h)$, 
for the source term $\S_h$ and also for the mesh velocity $\mathbf{V}_h$. 

Next, we have to update the mesh \textit{globally}. Let us denote with $\mathcal{V}_k$ the neighborhood of vertex number $k$, i.e. all those elements that have in common the node number $k$. The number of elements in 
the neighborhood  $\mathcal{V}_k$ is denoted with $N_k$. Since the velocity of each vertex is defined by the local predictor within each element, one has to deal with several, in general different, velocities for the 
same node, since all elements belonging to $\mathcal{V}_k$ will in general give a different velocity contribution, according to their element--local predictor. Since we do not admit the geometry to be discontinuous, 
we decide to fix a \textit{unique} node velocity $\overline{\mathbf{V}}_k^n$ to move the node. The final velocity is chosen to be the \textit{average velocity} considering all the contributions 
$\overline{\mathbf{V}}_{k,j}^n$ of the vertex neighborhood as  
\begin{equation}
\overline{\mathbf{V}}_k^n = \frac{1}{N_k}\sum \limits_{T_j^n \in \mathcal{V}_k}{\overline{\mathbf{V}}_{k,j}^n}, \qquad \textnormal{ with } \qquad 
\overline{\mathbf{V}}_{k,j}^n = \left( \int \limits_{0}^{1} \theta_l(\xi_{e,m(k)}, \eta_{e,m(k)}, \tau) d \tau \right) \widehat{\mathbf{V}}_{l,j}. 
\label{NodesVel}
\end{equation} 
The $\xi_{e,m(k)}$ and $\eta_{e,m(k)}$ are the vertex coordinates of the reference triangle $T_e$ corresponding to vertex number $k$, hence $m=m(k)$ with $1 \leq m \leq 3$ is a mapping from the global node number $k$ 
to the element--local vertex number. Since each node now has its own unique velocity, the vertex coordinates can be moved according to 
\begin{equation} 
	\mathbf{X}^{n+1}_{k}	= \mathbf{X}^{n}_{k}	+ \Delta t \, \overline{\mathbf{V}}_k^n,  
	\label{eqn.vertex.update}
\end{equation} 
and we can update all the other geometric quantities needed for the computation, e.g. normal vectors, volumes, side lengths, barycenter position, \textit{etc.}
\subsection{Finite Volume Scheme}
\label{sec.SolAlg}

The conservation law (\ref{PDE}) can be easily rewritten in a more compact space--time divergence form, which reads 
\begin{equation}
\tilde \nabla \cdot \tilde{\Q} = \mathbf{S}(\Q), 
\label{PDEdiv3D}
\end{equation} 
with
\begin{equation}
\tilde \nabla  = \left( \frac{\partial}{\partial x}, \, \frac{\partial}{\partial y}, \, \frac{\partial}{\partial t} \right)^T,  \qquad 
\tilde{\Q}  = \left( \mathbf{F}, \, \Q \right) = \left( \mathbf{f}, \, \mathbf{g}, \, \Q \right).
\end{equation}

Integration over the space--time control volume $C^n_i = T_i(t) \times \left[t^{n}; t^{n+1}\right]$ yields 
\begin{equation}
\int\limits_{t^{n}}^{t^{n+1}} \int \limits_{T_i(t)} \tilde \nabla \cdot \tilde{\Q} \, d\mathbf{x} dt = \int\limits_{t^{n}}^{t^{n+1}} \int \limits_{T_i(t)} \S \, d\mathbf{x} dt.   
\label{STPDE}
\end{equation} 
Applying Gauss' theorem allows us to write the left space--time volume integral above as sum of the flux  integrals computed over the space--time surface $\partial C^n_i$, hence 
\begin{equation}
\int \limits_{\partial C^n_i} \tilde{\Q} \cdot \ \mathbf{\tilde n} \, \, dS = 
\int\limits_{t^{n}}^{t^{n+1}} \int \limits_{T_i(t)} \S \, d\mathbf{x} dt,   
\label{I1}
\end{equation}    
where $\mathbf{\tilde n} = (\tilde n_x,\tilde n_y,\tilde n_t)$ is the outward pointing space--time unit normal vector on the space--time surface $\partial C^n_i$. 

\begin{figure}[!htbp]
	\begin{center} 
	\includegraphics[width=0.85\textwidth]{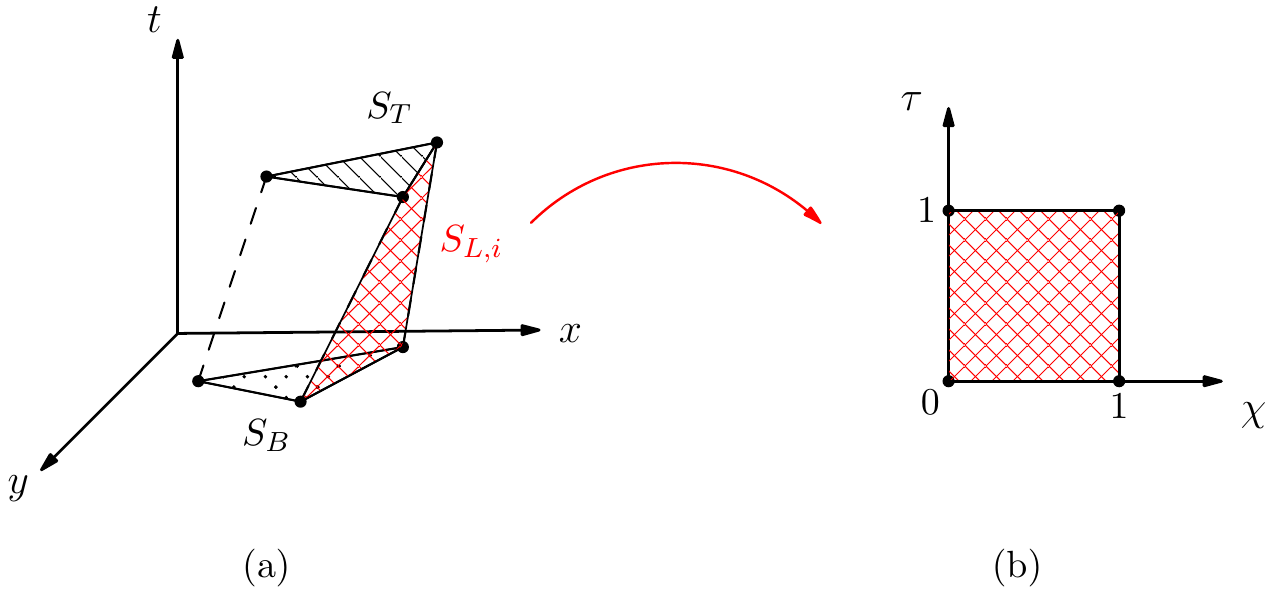}
	\caption{Physical space--time element (a) and parametrization of the lateral 
	         space--time subsurface $\partial C_{ij}^n$  (b).}
	\label{fig:STelem}
	\end{center} 
\end{figure}

The space--time surface $\partial C^n_i$ above involves overall five space--time sub--surfaces, as 
depicted in  Figure \ref{fig:STelem}:  
\begin{equation}
\partial C^n_i = \left( \bigcup \limits_{T_j(t) \in \mathcal{N}_i} \partial C^n_{ij} \right) 
\,\, \cup \,\, T_i^{n} \,\, \cup \,\, T_i^{n+1},  
\label{dCi}
\end{equation}    
where $\mathcal{N}_i$ denotes the so--called \textit{Neumann neighborhood} of triangle $T_i(t)$, i.e. 
the set of directly adjacent triangles $T_j(t)$ that share a common edge $\partial T_{ij}(t)$ with 
triangle $T_i(t)$. 
The common space--time edge $\partial C^n_{ij}$ during the time interval $[t^n;t^{n+1}]$ is denoted 
above by $\partial C^n_{ij} = \partial T_{ij}(t) \times [t^n;t^{n+1}]$.

The upper space--time sub--surface $T_i^{n+1}$ and the lower space--time sub--surface $T_i^{n}$ are 
parametrized by $0 \leq \xi \leq 1 \wedge 0 \leq \eta \leq 1-\xi$ and the mapping \eqref{xietaTransf}. 
They are orthogonal to the time coordinate, hence for these faces the space--time unit normal vectors 
simply read $\mathbf{\tilde n} = (0,0,1)$ for $T_i^{n+1}$ and $\mathbf{\tilde n} = (0,0,-1)$ for $T_i^{n}$, 
respectively. 

The lateral space--time sub--faces $\partial C^n_{ij}$ are defined using a simple bilinear 
parametrization, since the old vertex coordinates $\mathbf{X}_{ik}^n$ are given and the new ones 
$\mathbf{X}_{ik}^{n+1}$ are known from \eqref{eqn.vertex.update}.  
\begin{equation}
\partial C_{ij}^n = \mathbf{\tilde{x}} \left( \chi,\tau \right) = 
 \sum\limits_{k=1}^{4}{\beta_k(\chi,\tau) \, \mathbf{\tilde{X}}_{ij,k}^n },	
 \qquad 0 \leq \chi \leq 1,  \quad	0 \leq \tau \leq 1, 										 
\label{SurfPar}
\end{equation}
where $(\chi,\tau)$ represents a side-aligned local reference system according to Figure \ref{fig:STelem}. The $\mathbf{\tilde{X}}_{ij,k}^n$ are the  
physical space--time coordinate vectors for the four vertices that define the lateral space--time sub--surface  $\partial C_{ij}^n$. If $\mathbf{X}^n_{ij,1}$ and $\mathbf{X}^n_{ij,2}$ denote the two spatial nodes at time
$t^n$ that define the common spatial edge $\partial T_{ij}(t^n)$, then the four vectors  $\mathbf{\tilde{X}}_{ij,k}^n$ are  given by 
\begin{equation}
\mathbf{\tilde{X}}_{ij,1}^n = \left( \mathbf{X}^n_{ij,1}, t^n \right), \qquad 
\mathbf{\tilde{X}}_{ij,2}^n = \left( \mathbf{X}^n_{ij,2}, t^n \right), \qquad 
\mathbf{\tilde{X}}_{ij,3}^n = \left( \mathbf{X}^{n+1}_{ij,2}, t^{n+1} \right), \qquad 
\mathbf{\tilde{X}}_{ij,4}^n = \left( \mathbf{X}^{n+1}_{ij,1}, t^{n+1} \right).  
\label{eqn.lateralnodes} 
\end{equation} 
The $\beta_k(\chi,\tau)$ 
are a set of bilinear basis functions, which are defined as 
\begin{eqnarray}
 \beta_1(\chi,\tau) &=& (1-\chi)(1-\tau), \nonumber\\
 \beta_2(\chi,\tau) &=& \chi(1-\tau), \nonumber\\
 \beta_3(\chi,\tau) &=& \chi\tau, \nonumber\\
 \beta_4(\chi,\tau) &=& (1-\chi)\tau.
 \label{BetaBaseFunc}
\end{eqnarray}
From \eqref{eqn.lateralnodes} and \eqref{BetaBaseFunc} it follows that the temporal mapping is again simply
$t = t^n + \tau \, \Delta t$, hence $t_\chi = 0$ and $t_\tau = \Delta t$. 

The determinant of the coordinate transformation and the resulting space--time unit normal vector $\mathbf{\tilde n}_{ij}$ of  the sub--surface $\partial C_{ij}^n$ can be computed as follows: 
\begin{equation}
| \partial C_{ij}^n| = \left| \frac{\partial \mathbf{\tilde{x}}}{\partial \chi} \times \frac{\partial \mathbf{\tilde{x}}}{\partial \tau} \right|, 
\quad 
\mathbf{\tilde n}_{ij} = \left( \frac{\partial \mathbf{\tilde{x}}}{\partial \chi} \times \frac{\partial \mathbf{\tilde{x}}}{\partial \tau}\right) / | \partial C_{ij}^n|.
\label{n_lateral}
\end{equation}

Therefore, the ALE--type one--step finite volume scheme takes the following form: 
\begin{equation}
|T_i^{n+1}| \, \Q_i^{n+1} = |T_i^n| \, \Q_i^n - \sum \limits_{T_j \in \mathcal{N}_i} \,\, {\int \limits_0^1 \int \limits_0^1 
| \partial C_{ij}^n| \tilde{\Q}_{ij} \cdot \mathbf{\tilde n}_{ij} \, d\chi d\tau}
+ \int\limits_{t^{n}}^{t^{n+1}} \int \limits_{T_i(t)} \S(\mathbf{q}_h) \, d\mathbf{x} dt, 
\label{PDEfinal}
\end{equation}
where $\tilde{\Q}_{ij} \cdot \mathbf{\tilde n}_{ij}$ is an Arbitrary--Lagrangian--Eulerian numerical flux function to resolve the discontinuity of the predictor solution $\mathbf{q}_h$ at the space--time sub--face 
$\partial C_{ij}^n$. The surface integrals appearing in \eqref{PDEfinal} are approximated using multidimensional Gaussian quadrature rules, see \cite{stroud} for details. At the interface $\partial C_{ij}^n$  let us denote the local space--time predictor solution inside element $T_i(t)$ by $\q_h^-$ and 
the element--local predictor solution of the neighbor element $T_j(t)$ by $\q_h^+$, then a simple 
ALE--type Rusanov flux \cite{Dumbser2012} is given by 
\begin{equation}
  \tilde{\Q}_{ij} \cdot \mathbf{\tilde n}_{ij} =  
  \frac{1}{2} \left( \tilde{\Q}(\q_h^+) + \tilde{\Q}(\q_h^-)  \right) \cdot \mathbf{\tilde n}_{ij}  - 
  \frac{1}{2} s_{\max} \left( \q_h^+ - \q_h^- \right),  
  \label{eqn.rusanov} 
\end{equation} 
where $s_{\max}$ is the maximum eigenvalue of the ALE Jacobian matrix in spatial normal direction,  
\begin{equation} 
\mathbf{A}^{\!\! \mathbf{V}}_{\mathbf{n}}(\Q)=\partial ( \mathbf{F} \cdot \mathbf{n} ) / \partial \Q - 
(\mathbf{V} \cdot \mathbf{n}) \,  \mathbf{I}, \qquad 
\textnormal{ with } \qquad   
\mathbf{n} = \frac{(\tilde n_x, \tilde n_y)^T}{\sqrt{\tilde n_x^2 + \tilde n_y^2}}.  
\end{equation} 
Here, $\mathbf{V} \cdot \mathbf{n}$ is the local normal mesh velocity and $\mathbf{I}$ is the $\nu \times \nu$ identity matrix. It can be easily verified that
\begin{equation}
	\mathbf{V} = \frac{1}{\Delta t} \left( \begin{array}{c} x_\tau \\ y_\tau \end{array} \right), \quad  
  \mathbf{\tilde n}_{ij} = \left( \begin{array}{c} 
  \phantom{-} y_\chi \Delta t \\ 
  -x_\chi \Delta t \\    
   x_\chi y_\tau - y_\chi x_\tau 
   \end{array} \right), \quad \textnormal{ hence } \quad 
  \mathbf{V} \cdot \mathbf{n} = -\frac{\tilde n_t}{\sqrt{\tilde n_x^2 + \tilde n_y^2}}.   
\label{eqn.vnormal}    
\end{equation} 

A more sophisticated Osher--type ALE flux has been introduced in the Eulerian case in \cite{OsherUniversal} 
and has been employed in the Lagrangian framework in \cite{Dumbser2012}. It reads
\begin{equation}
  \tilde{\Q}_{ij} \cdot \mathbf{\tilde n}_{ij} =  
  \frac{1}{2} \left( \tilde{\Q}(\q_h^+) + \tilde{\Q}(\q_h^-)  \right) \cdot \mathbf{\tilde n}_{ij}  - 
  \frac{1}{2} \left( \int \limits_0^1 \left| \mathbf{A}^{\!\! \mathbf{V}}_{\mathbf{n}}(\boldsymbol{\Psi}(s)) \right| ds \right) \left( \q_h^+ - \q_h^- \right),  
  \label{eqn.osher} 
\end{equation} 
with the straight--line segment path connecting the left and right state as follows, 
\begin{equation}
\boldsymbol{\Psi}(s) = \q_h^- + s \left( \q_h^+ - \q_h^- \right), \qquad 0 \leq s \leq 1.  
\label{eqn.path} 
\end{equation} 
In \eqref{eqn.osher} above, the usual definition of the matrix absolute value operator applies, i.e. 
\begin{equation}
 |\mathbf{A}| = \mathbf{R} |\boldsymbol{\Lambda}| \mathbf{R}^{-1},  \qquad |\boldsymbol{\Lambda}| = \textnormal{diag}\left( |\lambda_1|, |\lambda_2|, ..., |\lambda_\nu| \right),  
\end{equation} 
with the right eigenvector matrix $\mathbf{R}$ and its inverse $\mathbf{R}^{-1}$. According to \cite{OsherUniversal} the path integral 
appearing in \eqref{eqn.osher} is approximated using Gaussian quadrature rules of sufficient accuracy.   

\section{Test Problems} 
\label{sec.test} 

In order to validate the two--dimensional ALE finite volume scheme presented in the previous section, we present some computational test problems in the following. 
We consider the Euler equations of compressible gas dynamics, which read: 
\begin{equation}
\label{eulerCons}
\Q_t + \f_x + \g_y = \S(x,y,t)
\end{equation}
with
\begin{equation}
\label{eulerTerms}
\Q = \left( \begin{array}{c} \rho \\ \rho u \\ \rho v \\ \rho E \end{array} \right), \quad \f = \left( \begin{array}{c} \rho u \\ \rho u^2 + p \\ \rho uv \\ u(\rho E + p) \end{array} \right), \quad \g = \left( \begin{array}{c} \rho v \\ \rho uv \\ \rho v^2 + p  \\ v(\rho E + p) \end{array} \right).  
\end{equation}
Here, $\rho$ denotes the fluid density , $\mathbf{v}=(u,v)$ is the fluid's velocity vector and $\rho E$ represents the total energy density. The vector of source terms is denoted by $\S$ and the fluid pressure by $p$, given in  terms of the conserved quantities by the equation of state of an ideal gas as 
\begin{equation}
\label{eqn.eos} 
p = (\gamma-1)\left(\rho E - \frac{1}{2} \rho (u^2+v^2) \right),  
\end{equation}
with the ratio of specific heats $\gamma$. 

For each of the following test cases we choose the local mesh velocity as the local fluid velocity according to Eqn. \eqref{NodesVel}, i.e.  
\begin{equation}
 \mathbf{V} = \mathbf{v}. 
\end{equation} 

\subsection{Numerical Convergence Results} 
\label{sec.conv.Rates}

The numerical convergence studies for the two--dimensional Lagrangian finite volume scheme are carried out considering a smooth convected isentropic vortex, see \cite{HuShuVortex1999}. 
The initial condition is given in terms of primitive variables and it consists in a linear superposition of a homogeneous background field and some perturbations $\delta$: 
\begin{equation}
\label{ShuVortIC}
(\rho, u, v, p) = (1+\delta \rho, 1+\delta u, 1+\delta v, 1+\delta p).
\end{equation}  
We set the vortex radius $r^2=(x-5)^2+(y-5)^2$, the vortex strength $\epsilon=5$ and the ratio of specific heats $\gamma=1.4$. The perturbation of entropy $S=\frac{p}{\rho^\gamma}$ is assumed 
to be zero, while the perturbations of temperature $T$ and velocity $\mathbf{v}$ are given by 
\begin{equation}
\label{ShuVortDelta}
\left(\begin{array}{c} \delta u \\ \delta v \end{array}\right) = \frac{\epsilon}{2\pi}e^{\frac{1-r^2}{2}} \left(\begin{array}{c} -(y-5) \\ \phantom{-}(x-5) \end{array}\right), \quad \delta S = 0, \quad \delta T = -\frac{(\gamma-1)\epsilon^2}{8\gamma\pi^2}e^{1-r^2}. 
\end{equation} 
From \eqref{ShuVortDelta} it follows that the perturbations for density and pressure are given by 
\begin{equation}
\label{rhopressDelta}
\delta \rho = (1+\delta T)^{\frac{1}{\gamma-1}}-1, \quad \delta p = (1+\delta T)^{\frac{\gamma}{\gamma-1}}-1. 
\end{equation} 

\begin{figure}[!htbp]
	\centering
		\includegraphics[width=0.95\textwidth]{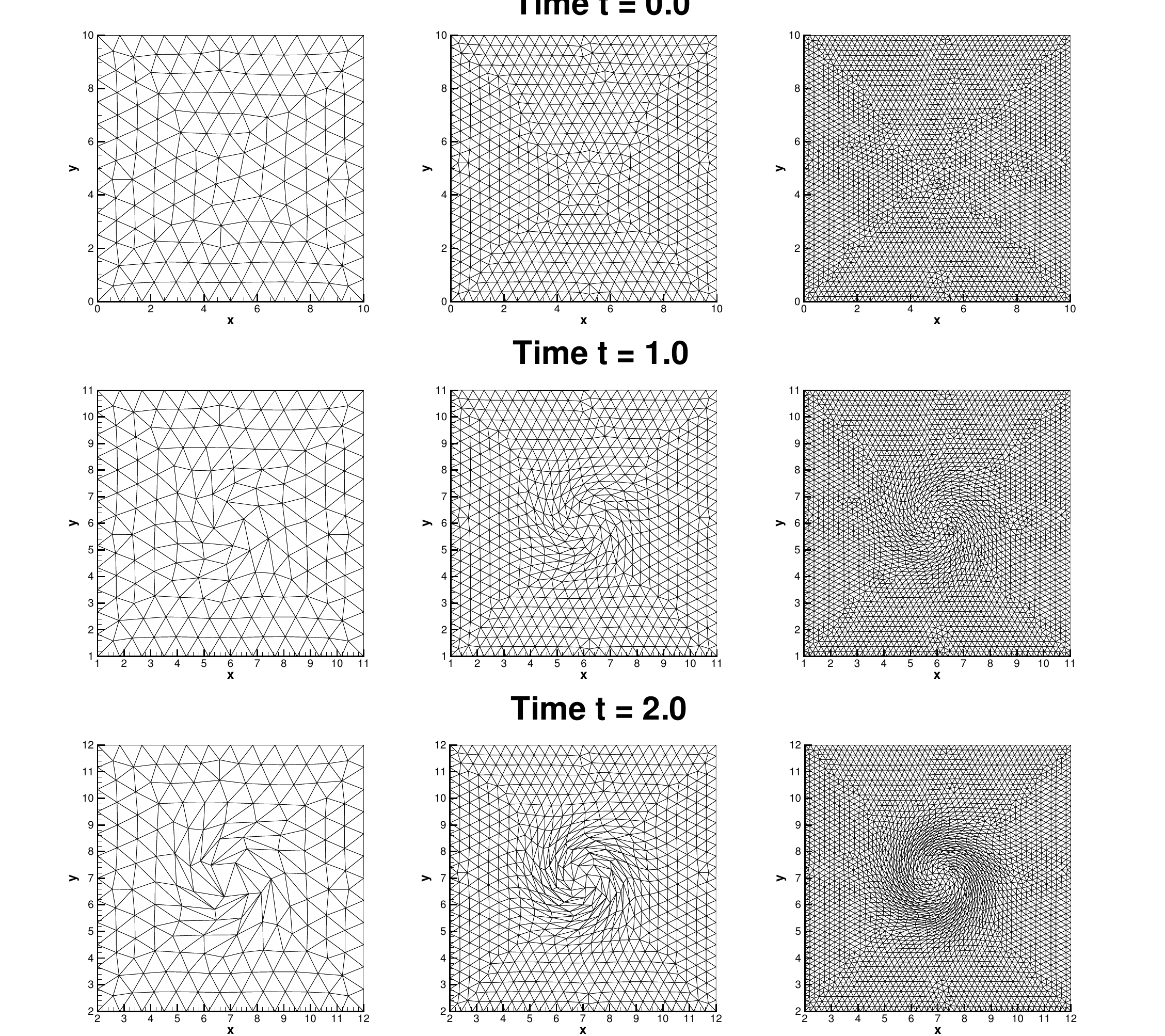}
		\caption{Different mesh sizes used for the convergence rates study at different time outputs: $t=0$ (top row), $t=1$ (middle row) and $t=2$ (bottom row). The total number of elements $N_E$ is increasing from the left grid ($N_E=320$), passing through the middle one ($N_E=1298$), to the right one ($N_E=5180$).}
	\label{fig:ConvRates}
\end{figure}

The initial computational domain $\Omega(0)=[0;10]\times[0;10]$ is square--shaped and it is surrounded by four periodic boundary conditions. The exact solution $\Q_e$ is the 
time--shifted initial condition given by $\Q_e(\x,t)=\Q(\x-\v_c t,0)$, with the convective mean velocity $\v_c=(1,1)$. We use the Osher--type numerical flux \eqref{eqn.osher} to run 
this test case until a final time of $t_f=1.0$ on successive refined meshes and for each mesh the corresponding error in $L_2$ norm is computed as 
\begin{equation}
  \epsilon_{L_2} = \sqrt{ \int \limits_{\Omega(t_f)} \left( \Q_e(x,y,t_f) - \w_h(x,y,t_f) \right)^2 dxdy },  
\end{equation} 
while the mesh size $h(\Omega(t_f))$ is taken to be the maximum diameter of the circumcircles of the triangles in the final domain $\Omega(t_f)$. The resulting numerical convergence 
rates are listed in Table \ref{tab.conv1}, whereas Figure \ref{fig:ConvRates} shows the temporal evolution of three different meshes used for the numerical convergence study. 

\begin{table}  
\caption{Numerical convergence results for the compressible Euler equations using the first up to sixth order version of the two--dimensional Lagrangian one--step WENO finite volume 
schemes presented in this article. The error norms refer to the variable $\rho$ (density) at time $t=1.0$.} 
\begin{center} 
\begin{small}
\renewcommand{\arraystretch}{1.0}
\begin{tabular}{ccccccccc} 
\hline
  $h(\Omega(t_f))$ & $\epsilon_{L_2}$ & $\mathcal{O}(L_2)$ & $h(\Omega,t_f)$ & $\epsilon_{L_2}$ & $\mathcal{O}(L_2)$ & $h(\Omega,t_f)$ & $\epsilon_{L_2}$ & $\mathcal{O}(L_2)$ \\ 
\hline
  \multicolumn{3}{c}{$\mathcal{O}1$} & \multicolumn{3}{c}{$\mathcal{O}2$}  & \multicolumn{3}{c}{$\mathcal{O}3$} \\
\hline
3.73E-01 & 9.525E-02 & -   & 3.43E-01  & 1.716E-02 & -   &  3.28E-01 & 1.614E-02 & -    \\ 
2.63E-01 & 6.907E-02 & 0.9 & 2.49E-01  & 1.109E-02 & 1.4 &  2.51E-01 & 6.943E-03 & 3.0  \\ 
2.14E-01 & 5.700E-02 & 0.9 & 1.69E-01  & 5.766E-03 & 1.7 &  1.68E-01 & 2.290E-03 & 2.7  \\ 
1.74E-01 & 4.752E-02 & 0.9 & 1.28E-01  & 3.027E-03 & 2.3 &  1.28E-01 & 9.274E-04 & 3.3  \\ 
\hline 
  \multicolumn{3}{c}{$\mathcal{O}4$} & \multicolumn{3}{c}{$\mathcal{O}5$}  & \multicolumn{3}{c}{$\mathcal{O}6$} \\
\hline
3.29E-01 & 4.717E-03 & -   & 3.29E-01  & 4.9463E-03 & -   & 3.29E-01 & 2.051E-03 & -    \\ 
2.51E-01 & 1.822E-03 & 3.5 & 2.51E-01  & 1.4648E-03 & 4.5 & 2.51E-01 & 5.803E-04 & 4.7  \\ 
1.67E-01 & 4.379E-04 & 3.5 & 1.67E-01  & 2.5937E-04 & 4.3 & 1.67E-01 & 8.317E-05 & 4.8  \\ 
1.28E-01 & 1.313E-04 & 4.4 & 1.28E-01  & 6.9664E-05 & 4.9 & 1.31E-01 & 1.994E-05 & 5.9  \\ 
\hline 
\end{tabular}
\end{small}
\end{center}
\label{tab.conv1}
\end{table} 

\subsection{Numerical Flux Comparison} 
\label{sec.flux comparison}
In order to study how the choice of the numerical flux does affect the solution, we consider the well-known Sod shock tube problem. The rectangular shaped computational domain $\Omega(0)=[-0.5;0.5]\times[-0.1;0.1]$ contains a total number of elements of $N_E=18018$. The initial condition reads
\begin{equation}
\Q(\x,0) = \left\{ \begin{array}{rlcc} \Q_L = &\left( 1.0,  0,0,1.0 \right) & \textnormal{ if } & \x \leq \x_0 \\ 
                                       \Q_R = &\left( 0.125,0,0,0.1 \right) & \textnormal{ if } & \x > \x_0        
                      \end{array}  \right. 
\label{eq:Sod_IC}
\end{equation}
with the position vector $\x_0=(0,y)$ which denotes the location of the discontinuity.

\begin{figure}[!htbp]
	\centering
		\includegraphics[width=0.7\textwidth]{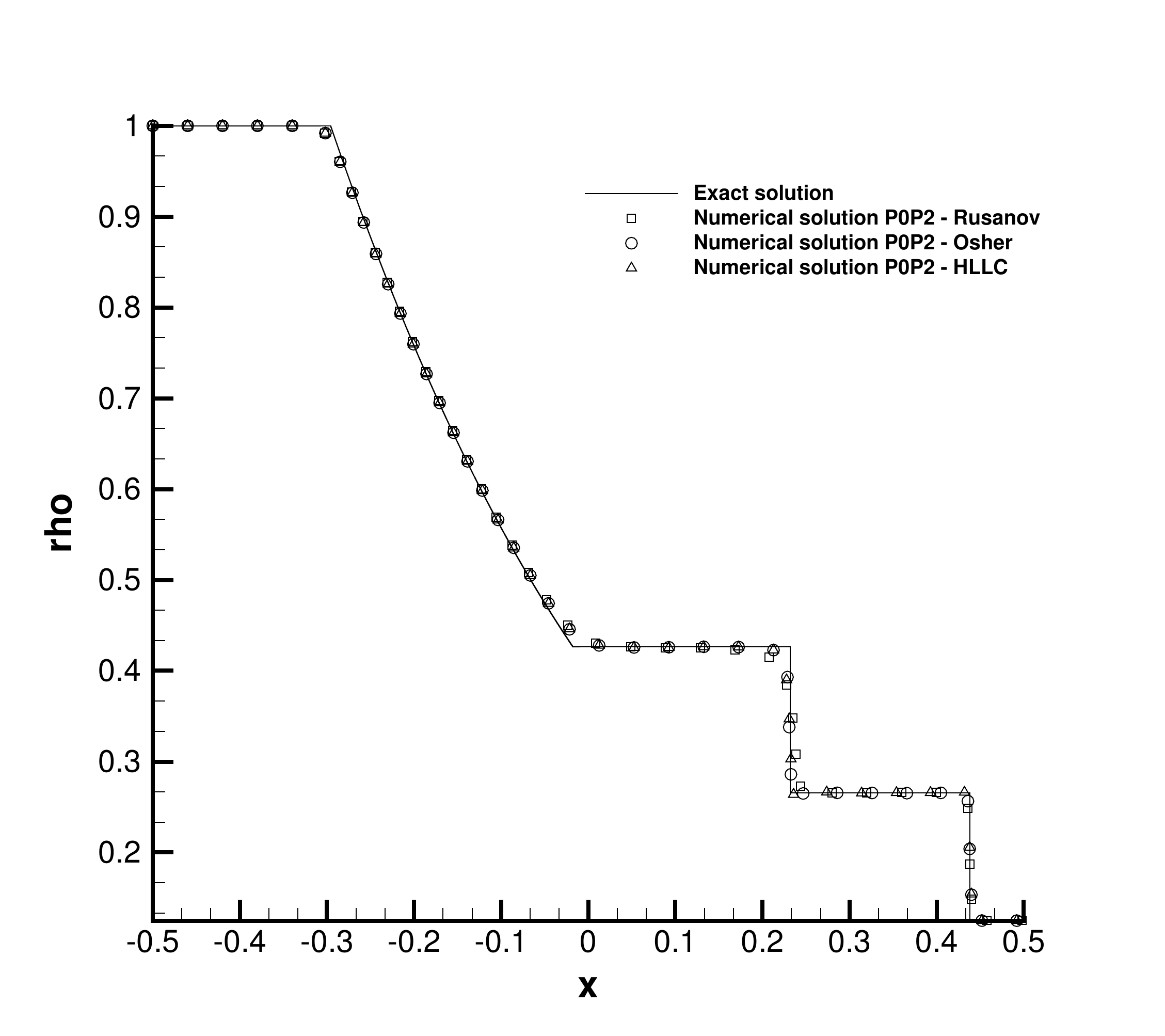}
		\caption{Sod shock tube test. Third order accurate numerical results obtained using three different numerical fluxes (Rusanov, Osher and HLLC) and comparison with the exact solution (solid line) at time $t=0.25$.}
	\label{fig:Sod}
\end{figure}

The Sod problem involves physically a rarefaction wave traveling towards left and both a contact and a shock wave which are moving to the right side of the domain. Exact solution is known from the usage of an exact Riemann solver, while the numerical results plotted in Figure \ref{fig:Sod} have been collected with three different numerical fluxes, namely the Rusanov--type and the Osher--type flux given by \eqref{eqn.osher} and \eqref{eqn.rusanov}, respectively. Furthermore we use also an HLLC--type numerical flux, firstly introduced for moving meshes by Van der Vegt et al. \cite{spacetimedg1,spacetimedg2} in the DG finite element framework. For a detailed description of the HLLC flux and its extension to dynamic grid motion we refer to \cite{spacetimedg1}.

Figure \ref{fig:Sod} shows a comparison between the exact solution and the numerical results at time $t=0.25$ obtained with a third order accurate scheme. The Rusanov flux is more diffusive if compared with the Osher and the HLLC fluxes, especially looking at the contact wave located at $x\approx0.23$, which tends to be smoothed, while on the contrary it is almost perfectly resolved by the Osher-- and the HLLC--ALE scheme.

Convergence rate studies for the three different numerical fluxes are reported in Table \ref{tab.conv2}. As test problem we use again the isentropic vortex fully described in section \ref{sec.conv.Rates}, but here we consider only a numerical scheme with an order of accuracy of $\mathcal{O}3$. The lowest $L_2-norm$ error is achieved by the very accurate Osher--type flux, as well as for the HLLC flux, whose error does not differ so much from the Osher scheme. The Rusanov type gives the highest error, as can be also observed in Figure \ref{fig:Sod}.

\begin{table}[!htbp]  
\caption{Comparison of numerical convergence results for the compressible Euler equations using the third order version of the two--dimensional Lagrangian one--step WENO finite volume schemes and three different types of numerical fluxes (Rusanov, Osher and HLLC). The error norms refer to the variable $\rho$ (density) at time $t=1.0$.}
\begin{center} 
\begin{small}
\renewcommand{\arraystretch}{1.0}
\begin{tabular}{ccccccccc} 
\hline
  \multicolumn{3}{c}{Rusanov} & \multicolumn{3}{c}{Osher}  & \multicolumn{3}{c}{HLLC} \\
\hline
  $h(\Omega(t_f))$ & $\epsilon_{L_2}$ & $\mathcal{O}(L_2)$ & $h(\Omega,t_f)$ & $\epsilon_{L_2}$ & $\mathcal{O}(L_2)$ & $h(\Omega,t_f)$ & $\epsilon_{L_2}$ & $\mathcal{O}(L_2)$ \\ 
\hline
3.61E-01 & 1.076E-01 & -   & 3.28E-01  & 1.614E-02 & -   &  3.31E-01 & 1.818E-02 & -    \\ 
2.51E-01 & 2.315E-02 & 4.2 & 2.51E-01  & 6.943E-03 & 3.2 &  2.51E-01 & 7.897E-03 & 3.0  \\ 
1.68E-01 & 8.658E-03 & 2.4 & 1.68E-01  & 2.290E-03 & 2.7 &  1.68E-01 & 2.621E-03 & 2.7  \\ 
1.28E-01 & 3.950E-03 & 2.9 & 1.28E-01  & 9.274E-04 & 3.3 &  1.28E-01 & 1.068E-03 & 3.3  \\ 
\hline 
\end{tabular}
\end{small}
\end{center}
\label{tab.conv2}
\end{table}

Looking at the computational time in Table \ref{tab.CPUtime}, the Rusanov flux is the most efficent scheme, but also the less accurate, while the Osher--type version gives the lowest error and is not the most expensive one, since the HLLC flux requires normally major computational efforts. We take Rusanov computational time as reference, hence $\eta$ evaluates simply the ratio between the current flux CPU time and the Rusanov one. Using HLLC one gets up to $\eta=1.7$, while at most $\eta=1.5$ is achieved with the Osher--type flux, which ensures the scheme to be the most accurate.    

\begin{table}[!htbp] 
\caption{Computational time for the convergence studies and the Sod shock tube results obtained with three different numerical fluxes (Rusanov, Osher and HLLC). CPU Time is measured in seconds $[s]$ and $\eta$ denotes the ratio w.r.t. the Rusanov computational time.}
\begin{center} 
\begin{small}
\renewcommand{\arraystretch}{1.0}
\begin{tabular}{ccccccc} 
\hline
                     & {Rusanov} &  &\multicolumn{2}{c}{Osher}  & \multicolumn{2}{c}{HLLC} \\
\hline
                     & CPU time  &  & CPU time & $\eta$         & CPU time & $\eta$  \\
\hline
                     & 6.21E+00 &   & 9.11E+00 & 1.5 & 1.06E+01  & 1.7  \\ 
\textit{Isentropic}  & 2.22E+01 &   & 2.49E+01 & 1.1 & 2.77E+01  & 1.2  \\ 
\textit{vortex}      & 9.64E+01 &   & 1.02E+02 & 1.1 & 1.12E+02  & 1.1  \\ 
                     & 1.55E+02 &   & 1.69E+02 & 1.1 & 1.93E+02  & 1.2  \\ 
\textit{Sod problem} & 1.22E+04 &   & 1.44E+04 & 1.2 & 1.59E+04  & 1.3  \\                     
\hline 
\end{tabular}
\end{small}
\end{center}
\label{tab.CPUtime}
\end{table}

\subsection{A Two--Dimensional Explosion Problem} 
\label{sec.RP2D}
This test problem is in some sense a two--dimensional extension of the classical one--dimensional shock tube problems. The initial domain $\Omega(0)$ is a circle of radius $R_o=1$ and the initial condition is given by 
\begin{equation}
  \Q(\x,0) = \left\{ \begin{array}{ccc} \Q_i & \textnormal{ if } & r \leq R, \\ 
                                        \Q_o & \textnormal{ if } & r > R,        
                      \end{array}  \right. 
\end{equation}
with $r^2=x^2+y^2$. The circle of radius $R=0.5$ separates the \textit{inner state} $\Q_i$ from the \textit{outer state} $\Q_o$. We assume $\gamma = 1.4$ and use the initial condition reported in Table \ref{tab:EP2D_IC}. 
 
\begin{table}[!htbp]
	\caption{Initial condition for the explosion test}
	\centering
		\begin{tabular}{ccccc}
		\hline
		            & $\rho$ & $u$ & $v$ & $p$ \\  
		\hline
		Inner state ($\Q_i$) & 1.0    & 0.0 & 0.0 & 1.0 \\ 
		Outer state ($\Q_o$) & 0.125  & 0.0 & 0.0 & 0.1 \\
		\hline
		\end{tabular}
	\label{tab:EP2D_IC}
\end{table}

In order to compute a reliable reference solution we assume rotational symmetry. Hence, one can simplify the multidimensional Euler equations to a one--dimensional system with geometric source terms, as proposed in \cite{ToroBook}, 
which reads 
\begin{equation}
\label{inhomEuler}
\Q_t + \F(\Q)_r = \S(\Q),
\end{equation}  
where
\begin{equation}
\label{matrixInhomEuler}
\Q = \left(\begin{array}{c} \rho \\ \rho u \\ \rho E \end{array}\right), \quad \F = \left(\begin{array}{c} \rho u \\ \rho u^2 + p \\  u(\rho E+p) \end{array}\right), \quad \S = -\frac{\alpha}{r}\left(\begin{array}{c} \rho u \\ \rho u^2 \\  u(\rho E+p) \end{array}\right).
\end{equation}
Here, $r$ denotes the radial direction and $u$ is the radial velocity, while $\alpha$ is a parameter that allows the system \eqref{inhomEuler} to be equivalent to the one--, two-- or three--dimensional Euler equations with rotational symmetry, according to its value: 
\begin{itemize}
  \item $\alpha =0$: plain 1D flow,
  \item $\alpha =1$: cylindrical symmetry (2D flow), 
  \item $\alpha =2$:	spherical symmetry (3D flow). 
\end{itemize}

\begin{figure}[!htbp]
	\centering
		\includegraphics[width=0.85\textwidth]{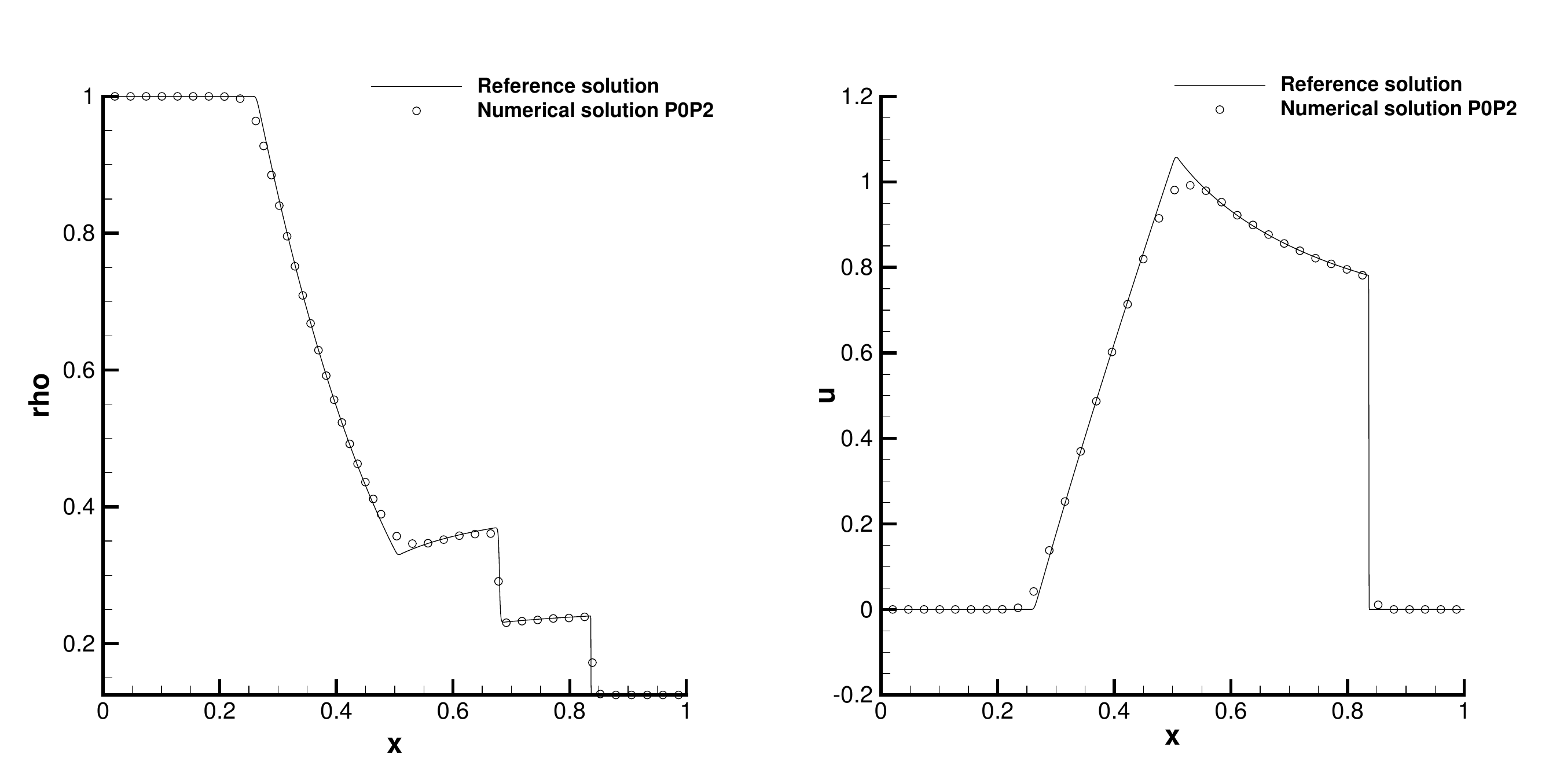}
	\caption{Profiles along the positive $x$-axis of density (left) and velocity (right) at time $t=0.2$. } 
	\label{fig:EP2DrhoU}
\end{figure}

\begin{figure}[!htbp]
	\centering
		\includegraphics[width=0.85\textwidth]{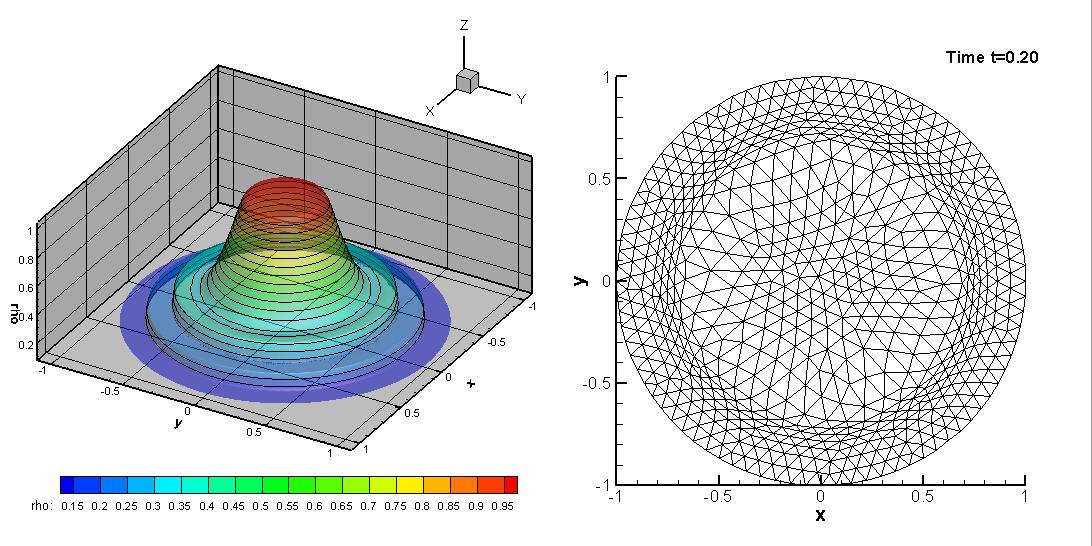}
	\caption{Left: density distribution at $t=0.20$ for the explosion test. Right: coarse distorted mesh at $t=0.20$.}
	\label{fig:EP2Drho3D}
\end{figure}

We choose $\alpha = 1$ for the 2D case and the inhomogeneous system of equations \eqref{inhomEuler} is solved using a second order MUSCL scheme with the Rusanov flux on a one--dimensional mesh of 15000 points in the radial 
interval $r \in [0;1]$. This solution is assumed to be our reference solution, which will be used to verify the accuracy of the two--dimensional Lagrangian finite volume scheme. Figure \ref{fig:EP2DrhoU} shows the comparison between the reference solution and the numerical solution obtained with the third order version of our Lagrangian one--step WENO scheme for density and velocity at time $t=0.20$: a circular shock wave is travelling away from the  center together with a contact wave, while a circular rarefaction wave is running towards the origin. The contact wave is very well resolved due to the use of the little diffusive Osher--type flux. The mesh used contained 68,324  triangles of characteristic mesh spacing $h=1/100$. A 3D view of the numerical solution as well as a very coarse version of the mesh are depicted in Fig. \ref{fig:EP2Drho3D}. 


\subsection{The Kidder Problem} 
\label{sec.Kidder}
The Kidder problem consists in an isentropic compression of a shell filled with an ideal gas. For this problem, an exact analytical solution has been proposed by Kidder in \cite{Kidder1976}. 
This is a very useful and widely used test \cite{Maire2009,Després2009}, to assure that a Lagrangian scheme does not produce spurious entropy. The shell has a time--dependent internal radius 
$r_i(t)$ and an external radius $r_e(t)$, while $r$ denotes the general radial coordinate with  $r_i(t) \leq r \leq r_e(t)$.  
The initial values for the internal and external radius are $r_i(0)=r_{i,0}=0.9$ and $r_e(0)=r_{e,0}=1.0$, respectively. 
In this test problem the ratio of specific heats is $\gamma=2$ and the initial density distribution is given by 
\begin{equation}
\rho_0 = \rho(r,0) = \left(\frac{r_{e,0}^2-r^2}{r_{e,0}^2-r_{i,0}^2}\rho_{i,0}^{\gamma-1}+\frac{r^2-r_{i,0}^2}{r_{e,0}^2-r_{e,0}^2}\rho_{e,0}^{\gamma-1}\right)^{\frac{1}{\gamma-1}} 
\end{equation}
with $\rho_{i,0}=1$ and $\rho_{e,0}=2$, which are the initial values of density defined at the internal and at the external frontier of the shell, respectively. 
The initial entropy $s_0= \frac{p_0}{\rho_0^\gamma} = 1$ is uniform, so that the initial pressure distribution is given by 
\begin{equation}
p_0(r) = s_0\rho_0(r)^\gamma. 
\end{equation}
Initially the fluid is at rest, hence $u=v=0$.

According to \cite{Kidder1976} we look for a solution of the the form $R(r,t) = h(t)r$, where $R(r,t)$ denotes the radius at time $t>0$ of a fluid particle initially located at radius $r$. Therefore the self-similar analytical solution for $t\in[0,\tau]$ reads
\begin{eqnarray}
\rho\left(R(r,t),t\right)&=&h(t)^{-\frac{2}{\gamma-1}}\rho_0\left[\frac{R(r,t)}{h(t)}\right], \nonumber\\
u_r\left(R(r,t),t\right)&=&\frac{d}{dt}h(t)\left[\frac{R(r,t)}{h(t)}\right], \nonumber\\
p\left(R(r,t),t\right)&=&h(t){-\frac{2\gamma}{\gamma-1}}p_0\left[\frac{R(r,t)}{h(t)}\right], 
\label{KidderAnalytical}
\end{eqnarray}
where the homothety rate is
\begin{equation}
h(t) = \sqrt{1-\frac{t^2}{\tau^2}}, 
\end{equation}
with $\tau$ representing the focalisation time 
\begin{equation}
\tau = \sqrt{\frac{\gamma-1}{2}\frac{(r_{e,0}^2-r_{i,0}^2)}{c_{e,0}^2-c_{i,0}^2}}, 
\end{equation}
and the internal and external sound speeds $c_i$ and $c_e$ are defined as 
\begin{equation}
c_i = \sqrt{\gamma\frac{p_i}{\rho_i}}, \qquad c_e = \sqrt{\gamma\frac{p_e}{\rho_e}}.
\end{equation}

The boundary conditions are imposed at the internal and the external frontier of the shell, where a space--time dependent state is prescribed according to the exact analytical solution \eqref{KidderAnalytical}. 

As done by Carr\'{e} et al. in \cite{Després2009} the final time of the simulation is chosen to be $t_f=\frac{\sqrt{3}}{2}\tau$, so that the resulting compression rate is $h(t_f)=0.5$. In this way the exact solution is given by a shell bounded with $0.45 \leq R \leq 0.5$.

\begin{figure}[!htbp]
	\centering
		\includegraphics[width=0.85\textwidth]{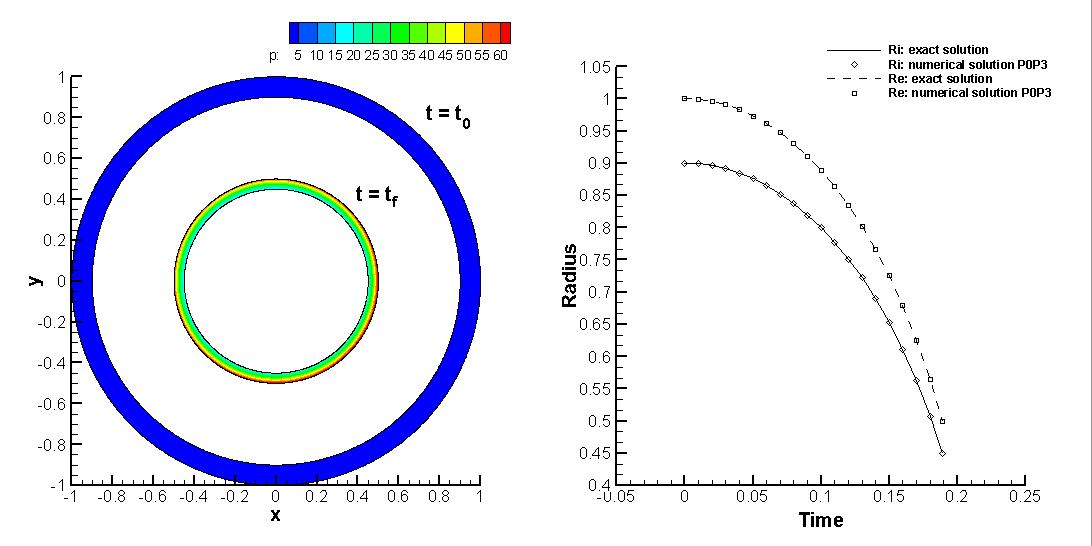}
	\caption{Left: Numerical solution at $t=0.4$ and at $t=t_f$. Right: Evolution of the internal and external radius of the shell and comparison between analytical and numerical solution.}
	\label{fig:Kidder}
\end{figure}

We use the Osher--type numerical flux \eqref{eqn.osher} to obtain the results shown in Figure \ref{fig:Kidder}, where we can observe an excellent agreement between the analytical and the numerical solution. 
The absolute error $|err|$ concerning the location of the internal and external radius at the end of the simulation has also be computed and is reported in Table \ref{tab:radiusKidder}. 

\begin{table}[!htbp]
	\begin{center}
		\begin{tabular}{|c|c|c|c|}
		\hline
		    		  			& $r_{ex}$ 		& $r_{num}$  & $|err|$     \\
		\hline
		\textit{Internal radius}	& 0.4500000 	& 0.4499936  & 6.40E-06 \\
		\hline
		\textit{External radius}	& 0.5000000 	& 0.4999922  & 7.80E-06 \\
		\hline
		\end{tabular}
	\end{center}
	\caption{Absolute error for the internal and external radius location between exact $R_{ex}$ and numerical $R_{num}$ solution.}
	\label{tab:radiusKidder}
\end{table}  

\subsection{The Saltzman Problem} 
\label{sec.Saltzman}
\begin{figure}[!htbp]
	\centering
		\includegraphics[width=0.75\textwidth]{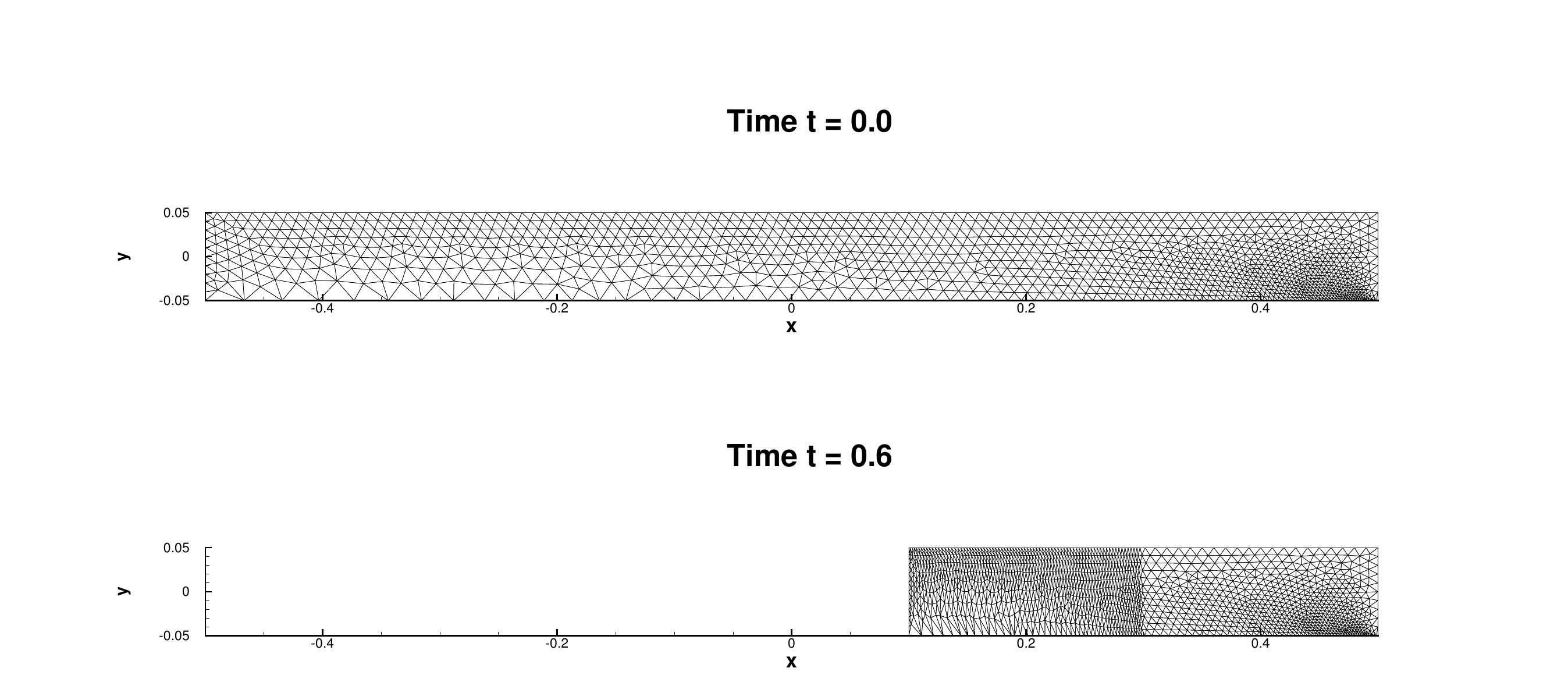}
	\caption{The initial and the final mesh for the Saltzman test problem.}
	\label{fig:SaltzMesh}
\end{figure}

A strong one--dimensional shock wave driven by a piston constitutes this test problem. As proposed in \cite{Maire2009,chengshu2}, the initial two--dimensional domain is a box $\Omega(0)=[0;1]\times[0;0.1]$, filled with an 
ideal gas. The piston is pushing the gas from the left to the right and it is moving with constant velocity $\mathbf{v}_p = (1,0)$. The computational domain is covered by a distorted unstructured mesh, composed of 
$100\times 10$ elements along the boundaries, as reported in Figure \ref{fig:SaltzMesh}.

The initial condition is given by an ideal gas at rest ($\mathbf{v}_0=\mathbf{0}$) with a ratio of specific heats $\gamma = \frac{5}{3}$, an initial density of $\rho_0=0$ and an internal energy of $e_0=10^{-6}$, 
corresponding to an initial pressure of $p_0=6.67 \cdot 10^{-7}$. 

As boundary conditions we set fixed slip wall boundaries on the upper and lower boundary, outflow on the right and a moving slip wall boundary condition for the piston on the left. According to \cite{chengshu2} we use 
initially a Courant number of $CFL=0.01$, in order to respect the geometric $CFL$ condition and to avoid mesh elements crossing over. After this initial phase, the CFL condition is raised to its usual value of $CFL=0.5$ 
in two space dimensions. 

The exact solution is a one--dimensional infinite strength shock wave and it can be computed by solving the Riemann problem given in Table \ref{tab:SaltzEx}. 
\begin{table}[!htbp]
	\centering
		\begin{tabular}{|c|c|c|}
		\hline
						& \textit{Left state} & \textit{Right state} \\
	  \hline
	  \hline
	  $\rho$ 	& 1.0  								& 1.0      				\\
	  $u$			& 1.0  								& -1.0    				\\
	  $v$			& 0.0  								& 0.0      				\\ 
	  $p$			& $6.67\cdot10^{-7}$	& $6.67\cdot10^{-7}$ \\
	  \hline
		\end{tabular}
	\caption{One--dimensional Riemann problem for obtaining the exact solution of the Saltzman problem.}
	\label{tab:SaltzEx}
\end{table}

The details of the algorithm that computes the exact solution of the Riemann problem are given in the book of Toro \cite{ToroBook}. The exact solution has then to be shifted by a certain quantity $d$ to the right, 
corresponding to the movement of the piston during the time of the simulation $t_f$, i.e. 
\begin{equation}
d = u_p\cdot t_f.
\end{equation} 
The exact post shock density is $\rho_e = 4.0$ and the final time is chosen to be $t_f = 0.6$, according to \cite{chengshu2}. Therefore, the exact final shock location is at $x=0.3$.  

\begin{figure}[!htbp]
	\centering
		\includegraphics[width=0.85\textwidth]{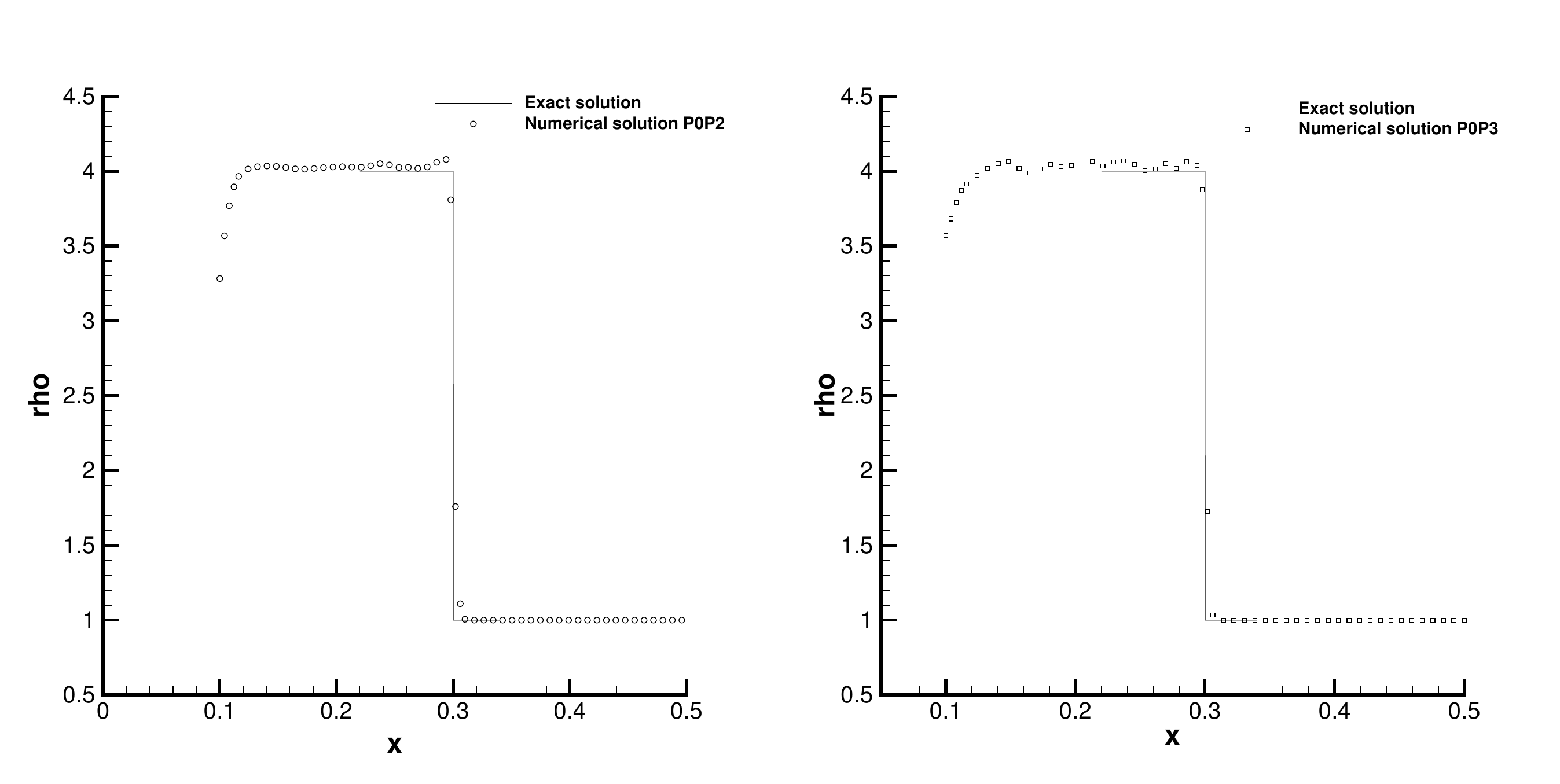}
		\caption{Numerical results for the density distribution of the Saltzman problem: third order scheme (left) and fourth order scheme (right). Both results are compared with the analytical solution, obtained by solving the Riemann problem described in Table \ref{tab:SaltzEx}.}
		\label{fig:Saltz_rhoEx.pdf}
\end{figure}

\begin{figure}[!htbp]
	\centering
		\includegraphics[width=0.85\textwidth]{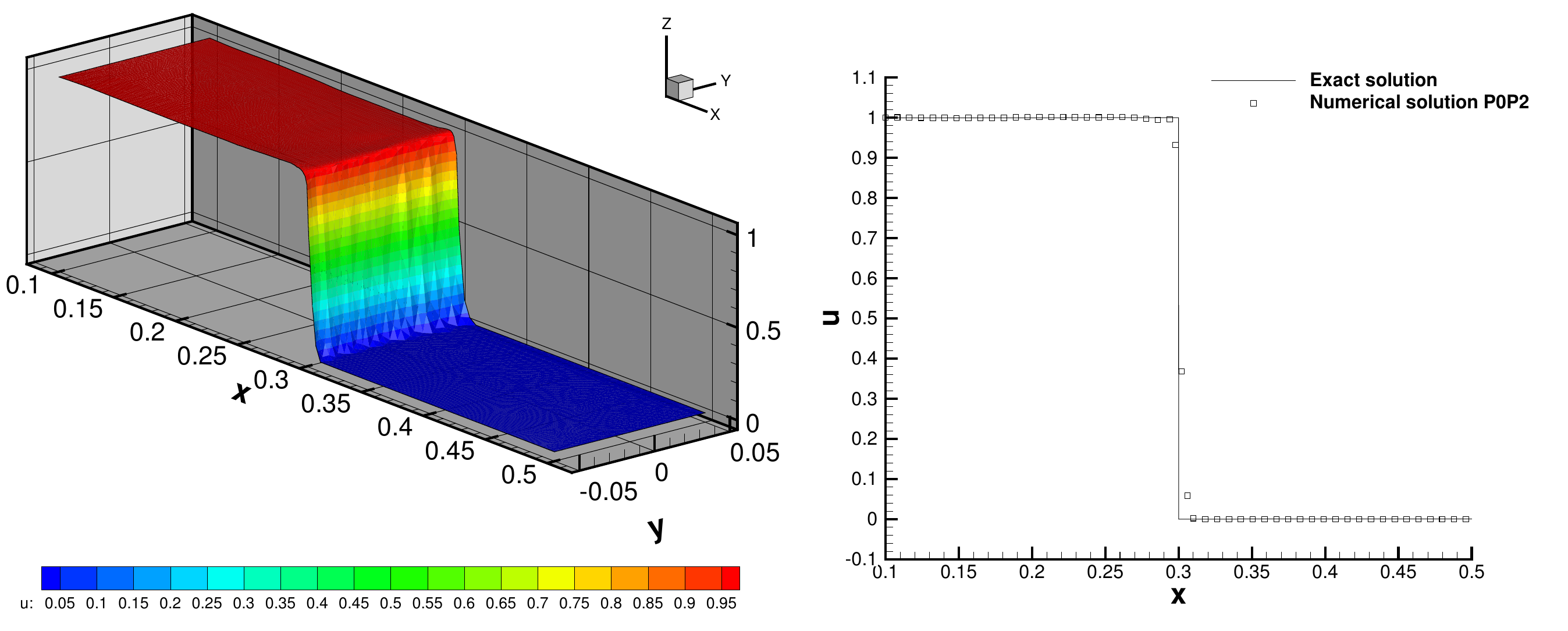}
		\caption{3D visualization for the velocity distribution obtained by a third order finite volume scheme (left) and comparison with the analytical solution (right).}
		\label{fig:SaltzP0P2_u3DEx.pdf}
\end{figure}

\begin{figure}[!htbp]
	\centering
		\includegraphics[width=0.85\textwidth]{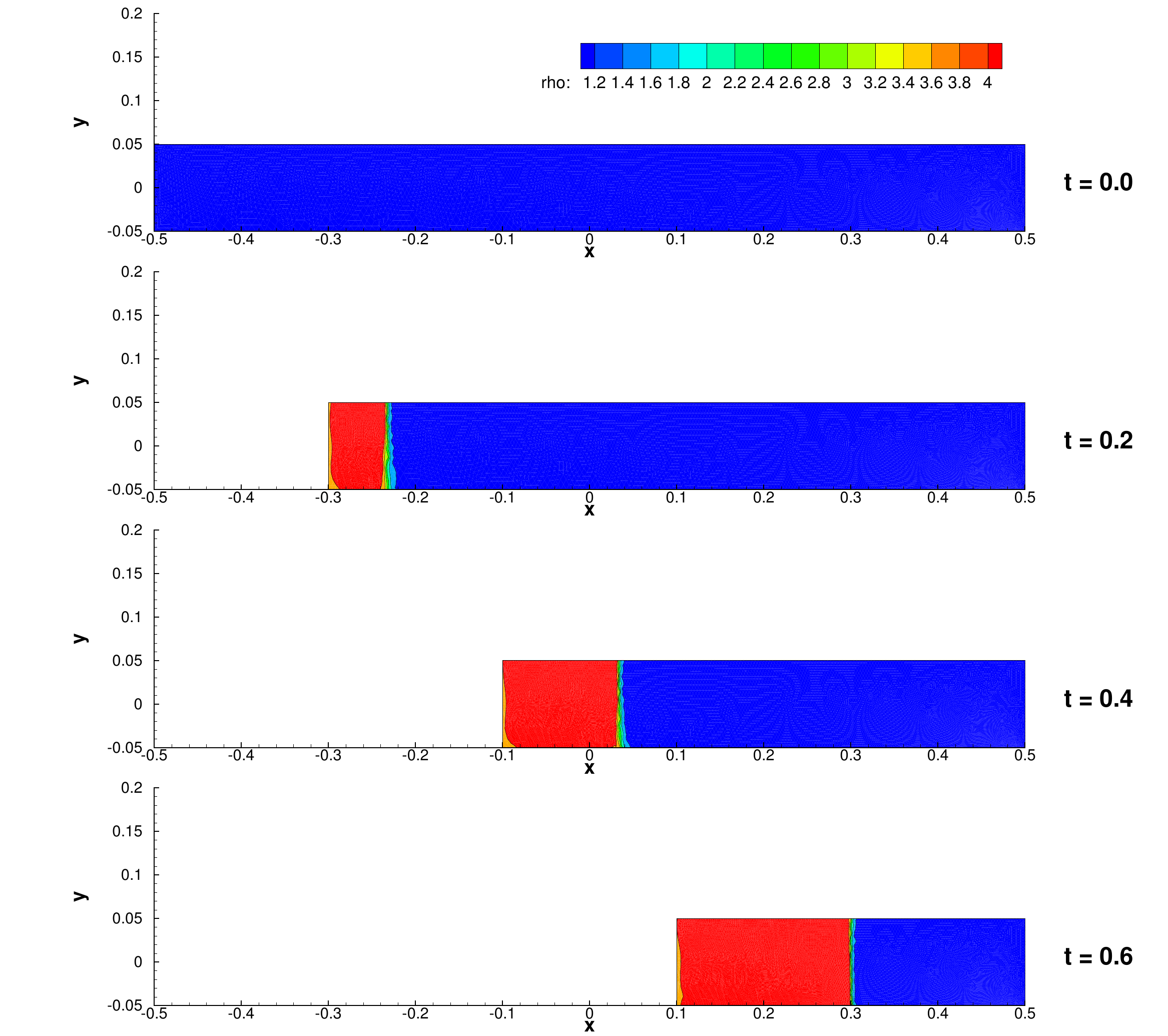}
		\caption{Progressive evolution of the numerical solution of the Saltzman problem (density distribution): the mesh is compressed by the piston on the left, while a shock wave is travelling towards the right. The final location  
		         coincides perfectly well with the exact solution.}
		\label{fig:SaltzP0P3_progressive.pdf}
\end{figure}

We use a Rusanov--type flux \eqref{eqn.rusanov} and the solution is computed, both, for a third order and a fourth order scheme. We can notice a very good agreement regarding the final density distribution, as depicted in Figure \ref{fig:Saltz_rhoEx.pdf}, and the solution for the velocity, where we observe a very sharp discontinuity at the shock location, as clearly shown in Figure \ref{fig:SaltzP0P2_u3DEx.pdf}. The decrease of the density near 
the piston in Figure \ref{fig:Saltz_rhoEx.pdf} is due to the well known \textit{wall--heating problem}, see \cite{toro.anomalies.2002}. One can notice from Figure \ref{fig:SaltzP0P3_progressive.pdf} the progressive mesh 
compression at the inflow, where the piston is pushing, and the shock wave that is travelling through the domain, reaching its final location exactly according to the reference solution ($x=0.3$).  

\subsection{Single Mach Reflection Problem}  
\label{sec.smr} 

\begin{figure}[!htbp]
\begin{center}
\includegraphics[width=0.85\textwidth]{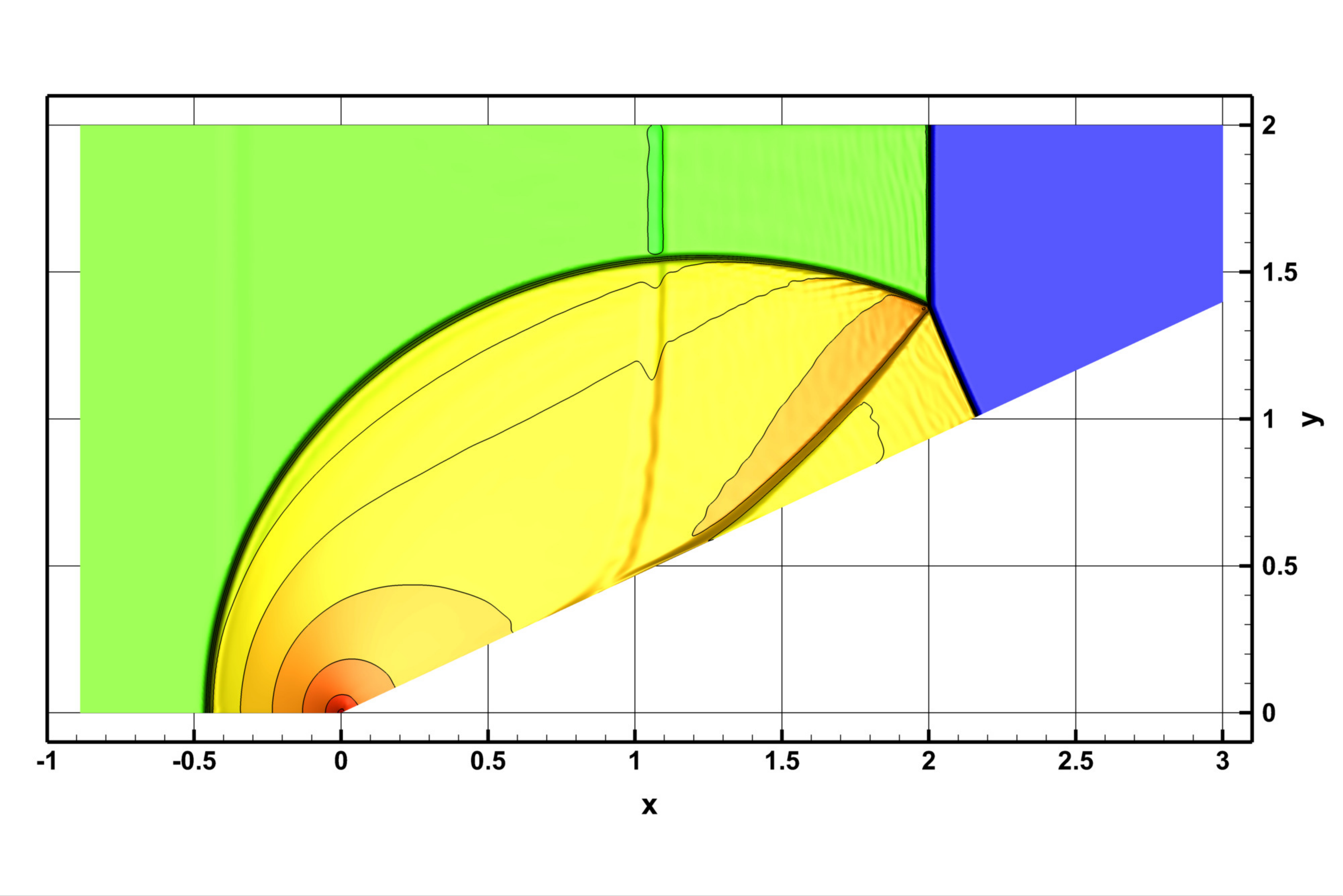}         
\includegraphics[width=0.75\textwidth]{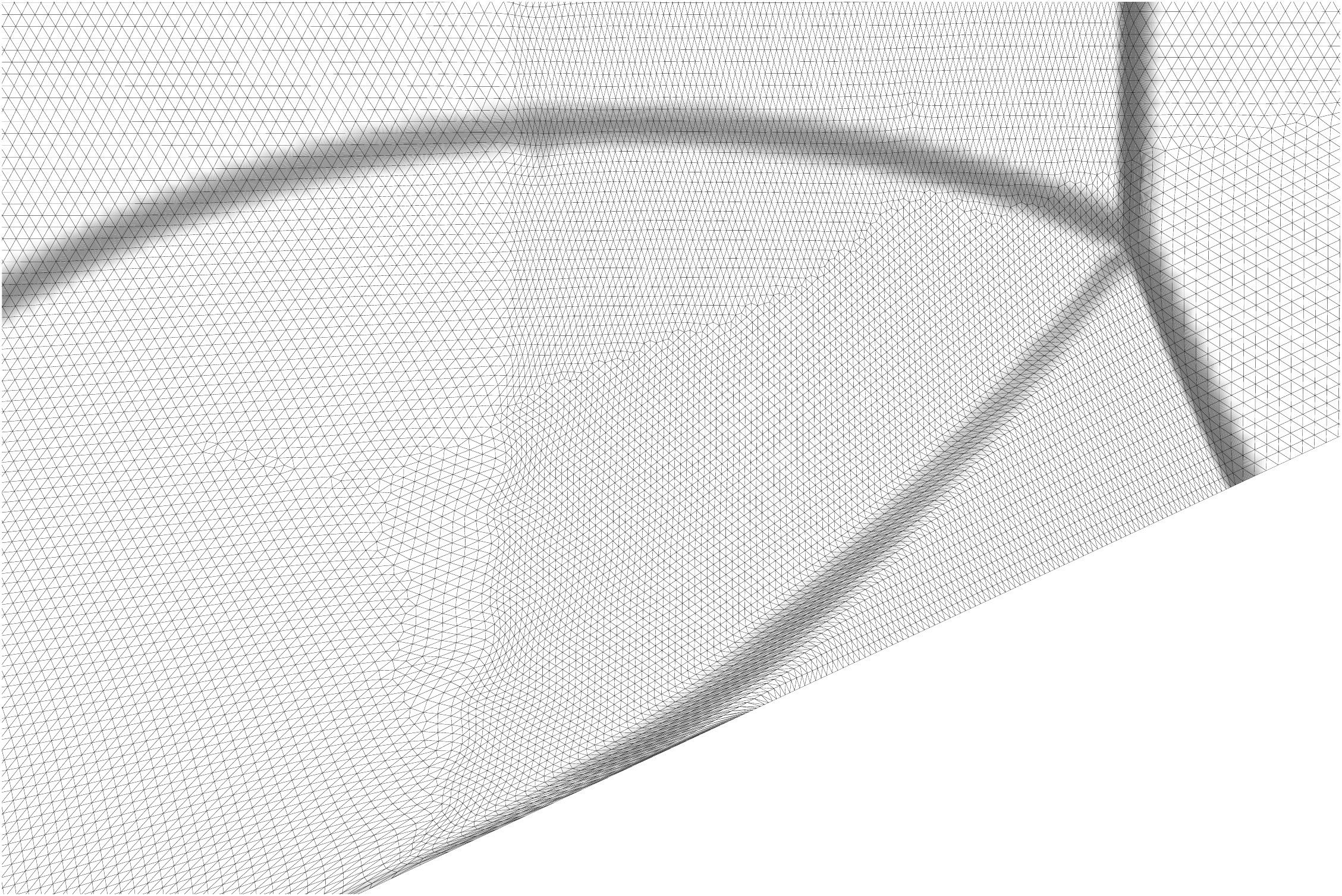}    
\caption{Top: 21 equidistant density contours from $\rho=1.2$ to $\rho=3$ for the single Mach reflection problem 
($M_s=1.7, \alpha=25^\circ$). Bottom: zoom into the mesh at the final time $t=1.2$.} 
\label{fig.smr}
\end{center}
\end{figure}

In this section we solve the single Mach reflection problem. It consists of a shock wave traveling
at a shock Mach number of $M_s=1.7$ that hits a wedge of an angle $\alpha=25^\circ$. Experimental reference data 
in form of Schlieren images as well as numerical reference solutions are shown, e.g. in \cite{ToroBook}. 
The upstream density and pressure are $\rho_0=1$ and $p_0 = 1/\gamma$, respectively, where the ratio
of specific heats is $\gamma = 1.4$. In the following, the indices $_0$ and $_1$ will denote upstream 
and downstream states, respectively. Using the Rankine-Hugoniot conditions, we setup the computation 
with a shock wave initially centered at $x=-0.04$ that travels at $M_s=1.7$ to the right into a medium at 
rest. The wedge is defined by the triangle $\mathcal{W}_{25}$, composed of the vertices $(0,0)$, 
$(3.0, 0)$ and $(3.0, 1.398923)$. The initial computational domain 
$\Omega(0) = [-2;3] \times [0;2] \, \backslash \, \mathcal{W}_{25}$ is discretized with 
177,440 triangles of characteristic size $h=1/100$ and a $P_0P_2$ WENO finite volume scheme is used. 
At the final time $t=1.2$ the exact shock location must be $x=2$. The upper and lower boundaries of 
the computational domain are discretized by slip walls. At the other boundaries Dirichlet boundary conditions 
consistent with the initial condition are imposed.  \\
The results of the computation are depicted in Fig. \ref{fig.smr}, where the density contour lines are shown, 
as well as a zoom into the final mesh, which is very distorted along the slip line. This is to verify that our 
proposed high order unstructured finite volume algorithm can handle reasonably strong mesh deformations, which 
are a common feature of Lagrangian computations. For even stronger mesh deformations, as they occur in the double
Mach reflection problem, a proper remeshing and projection strategy will be implemented in the future. The 
general flow features depicted in Fig. \ref{fig.smr} agree very well with the computations and the experimental data shown 
in \cite{ToroBook}. 

\subsection{Two-dimensional Riemann Problems}  
\label{sec.rp2d} 

\begin{figure}[!htbp]
\begin{center}
\begin{tabular}{cc} 
\includegraphics[width=0.47\textwidth]{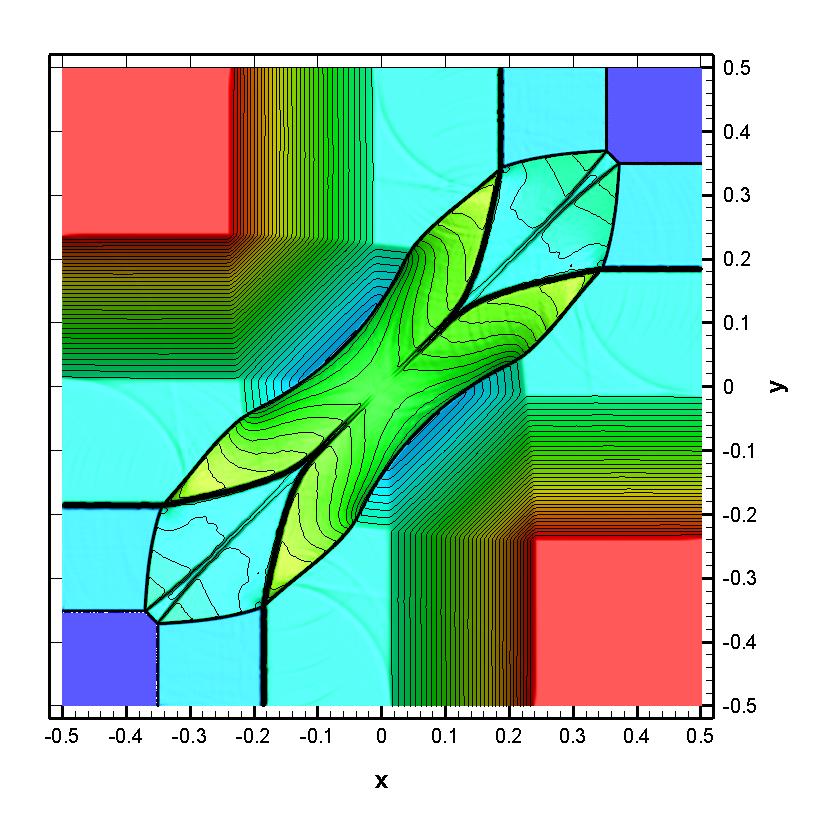}  &           
\includegraphics[width=0.47\textwidth]{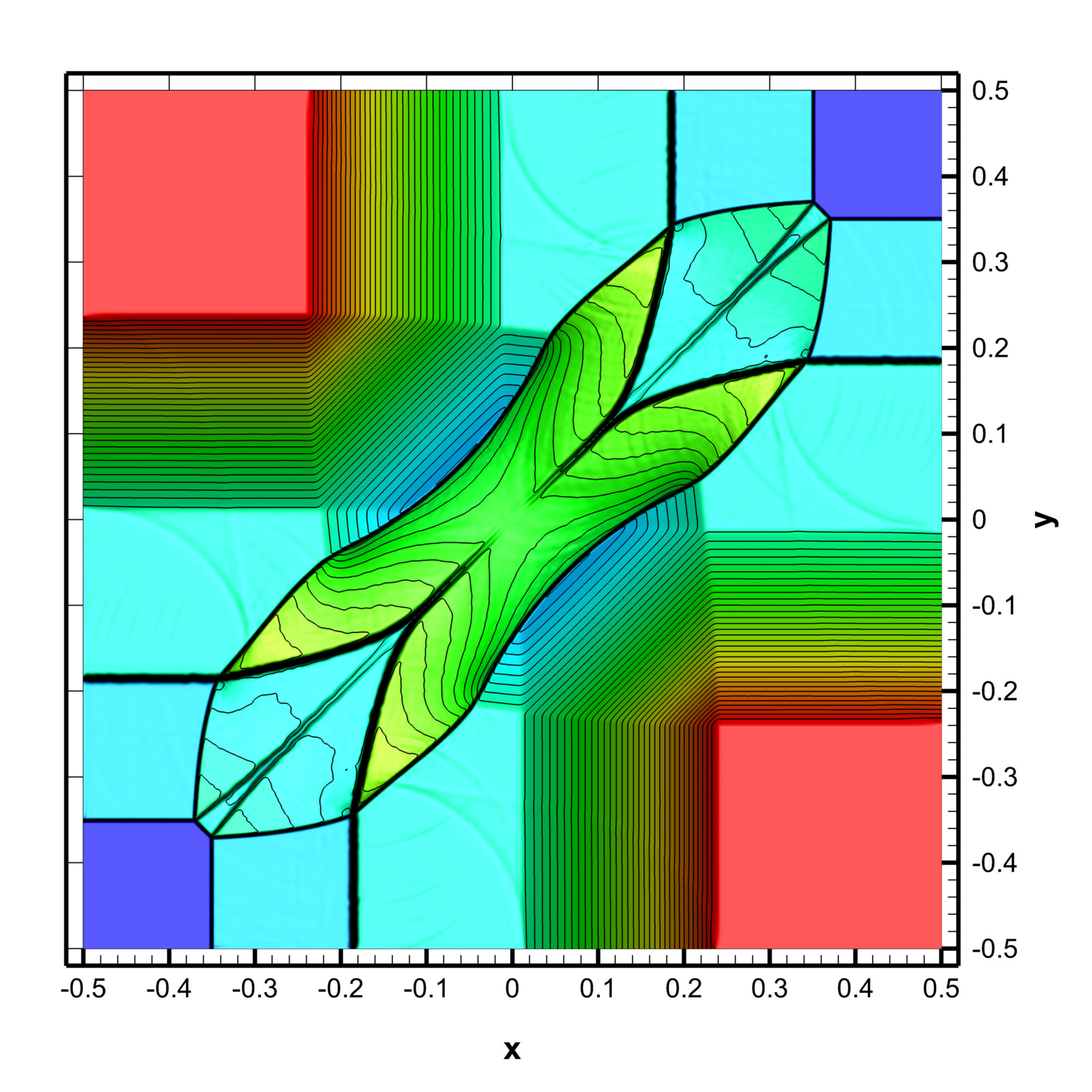} \\   
\multicolumn{2}{c}{
\includegraphics[width=0.70\textwidth]{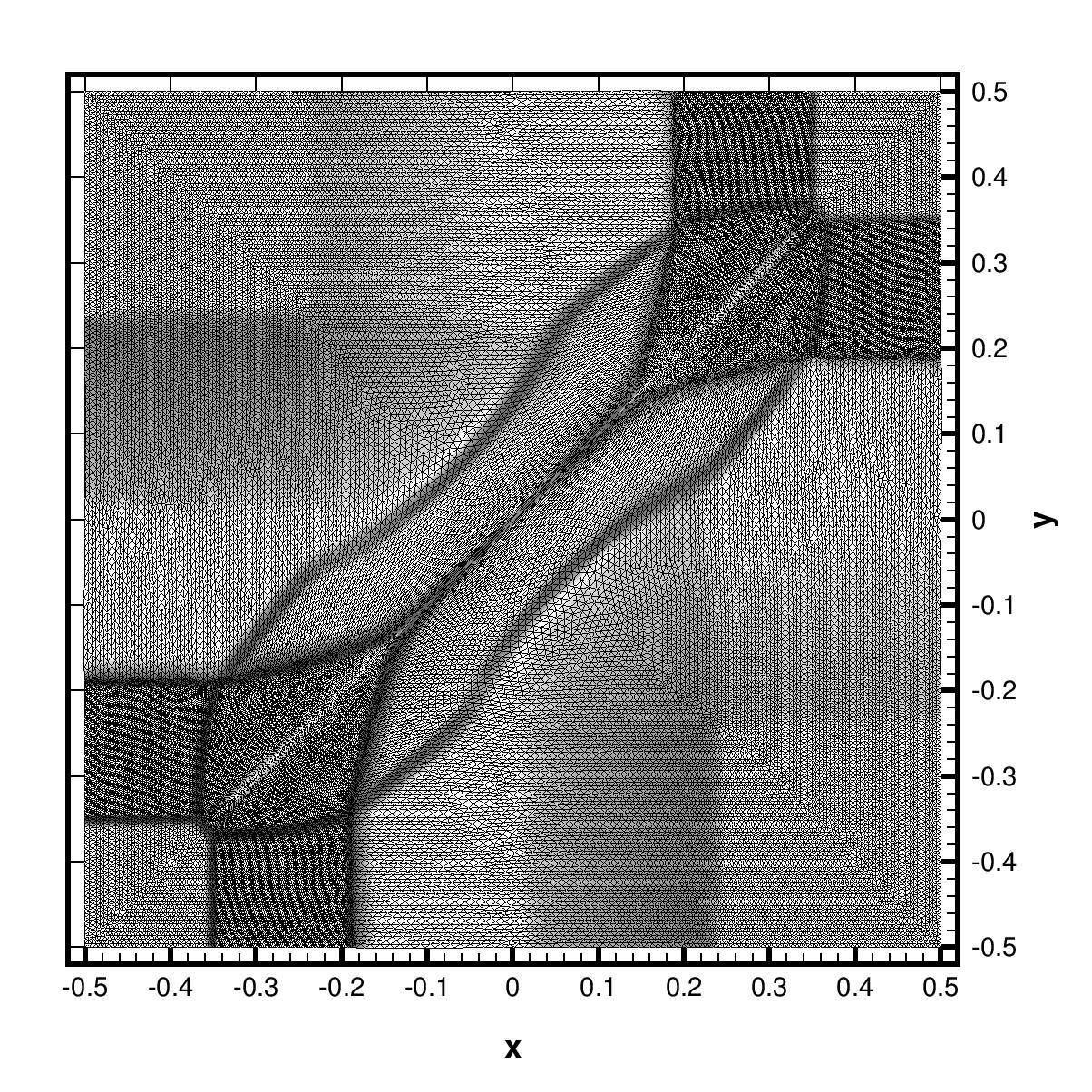}    
} 
\end{tabular} 
\caption{Top left: 50 equidistant density contour lines from $\rho=0.14$ to $\rho=0.98$ for the numerical solution obtained with our third order unstructured Lagrangian one--step WENO finite volume scheme 
for the two--dimensional Riemann problem RP1 at $t=0.2$. Top right: Same contour lines for the reference solution. Bottom: computational mesh of the Lagrangian scheme at the final time $t=0.2$.} 
\label{fig.rp1}
\end{center}
\end{figure}

\begin{figure}[!htbp]
\begin{center}
\begin{tabular}{cc} 
\includegraphics[width=0.47\textwidth]{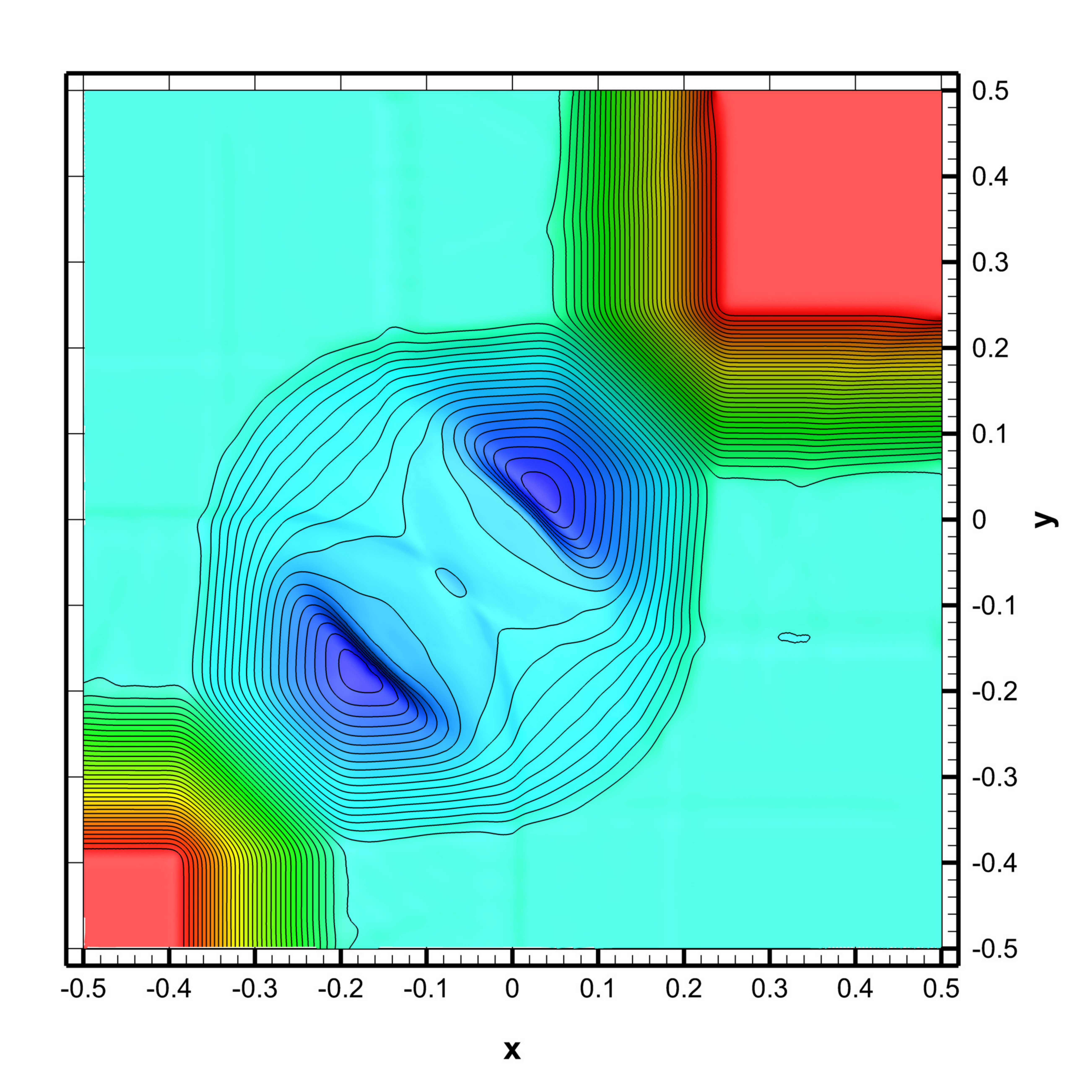}  &          
\includegraphics[width=0.47\textwidth]{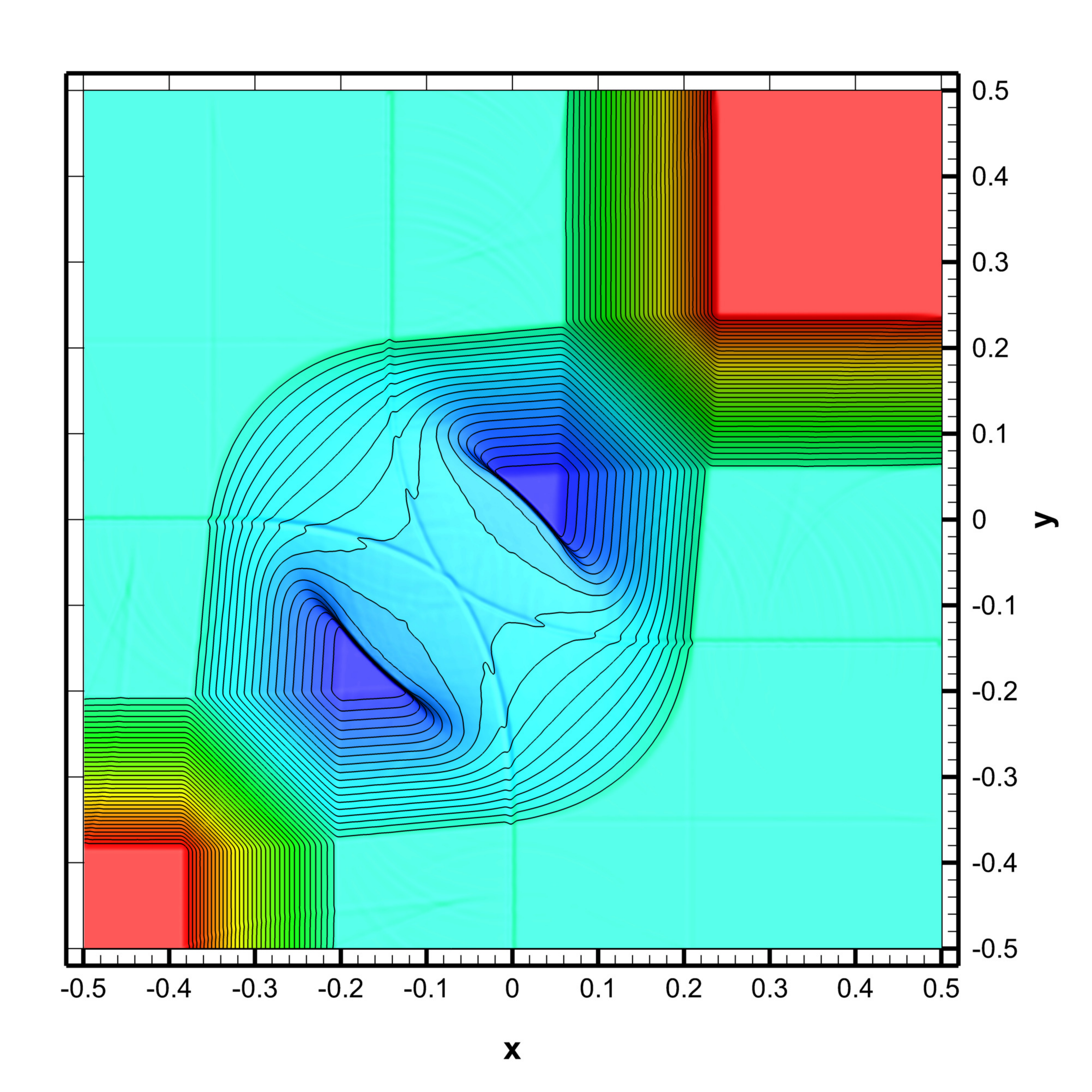} \\   
\multicolumn{2}{c}{
\includegraphics[width=0.70\textwidth]{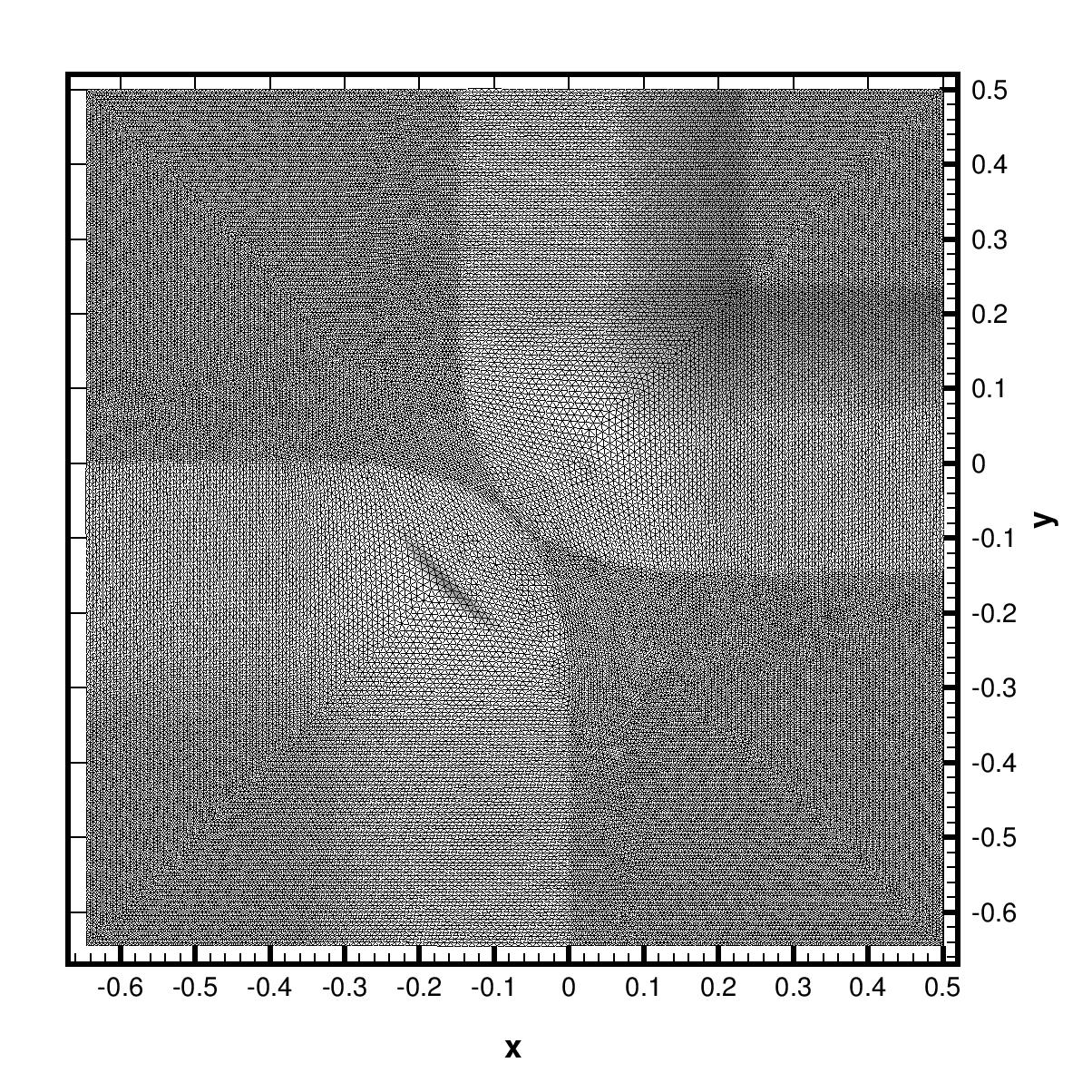}    
}
\end{tabular} 
\caption{Top left: 51 equidistant density contour lines from $\rho=0.27$ to $\rho=0.98$ for the numerical solution obtained with our third order unstructured Lagrangian one--step WENO finite volume scheme 
for the two--dimensional Riemann problem RP2 at $t=0.2$. Top right: Same contour lines for the reference solution. Bottom: computational mesh of the Lagrangian scheme at the final time $t=0.2$.} 
\label{fig.rp2}
\end{center}
\end{figure}

In this section we solve a set of two--dimensional Riemann problems, whose initial conditions are given by 
\begin{equation}
 \Q(x,0) = \left\{ \begin{array}{ccc} 
 \Q_1 & \textnormal{ if } & x > 0 \wedge y > 0,    \\ 
 \Q_2 & \textnormal{ if } & x \leq 0 \wedge y > 0, \\ 
 \Q_3 & \textnormal{ if } & x \leq 0 \wedge y \leq 0, \\ 
 \Q_4 & \textnormal{ if } & x > 0 \wedge y \leq 0.    
  \end{array} \right. 
\end{equation} 
A large set of such two--dimensional Riemann problems has been presented in great detail in the paper by Kurganov and Tadmor 
\cite{KurganovTadmor2002}. The initial conditions for the three configurations presented in this article are listed in 
Table \ref{tab.rp2d.ic}. The initial computational domain is defined as $\Omega(0)=[-0.5;0.5] \times [-0.5;0.5]$. 
The Lagrangian simulations are carried out with a third order one--step WENO scheme using an unstructured triangular mesh
composed of 90,080 elements with an initial characteristic mesh spacing of $h=1/200$. The reference solution is computed with 
a high order Eulerian one--step WENO finite volume scheme as presented in \cite{Dumbser20088209,Dumbser2007204,Dumbser2007693}, 
using a very fine mesh composed of 2,277,668 triangles with characteristic mesh spacing $h=1/1000$. In all cases, the Rusanov flux
has been used. The exact solutions of the one-dimensional Riemann problems are imposed as boundary conditions on the four boundaries 
of the domain. The results obtained with the Lagrangian-type scheme together with the final mesh and the Eulerian reference solution 
are depicted in Figures \ref{fig.rp1} - \ref{fig.rp2}. We can note a very good qualitative agreement of the Lagrangian solution with 
the reference solution, as well as with the results published in \cite{KurganovTadmor2002}. 

\begin{table}[!htbp]   
\caption{Initial conditions for the two--dimensional Riemann problems.} 
\begin{center} 
\renewcommand{\arraystretch}{1.0}
\begin{tabular}{c|cccc|cccc} 
\hline
\multicolumn{9}{c}{Problem RP1 (Quadruple Sod problem)} \\
\hline
     & \multicolumn{4}{c|}{$x \leq 0$} & \multicolumn{4}{|c}{$x>0$} \\
\hline
     & $\rho$ & $u$ & $v$  & p & $\rho$ & $u$ & $v$  & p  \\ 
\hline
$y > 0$    & 1.0   & 0.0 & 0.0 & 1.0 & 0.125 & 0.0 & 0.0 & 0.1 \\ 
$y \leq 0$ & 0.125 & 0.0 & 0.0 & 0.1 & 1.0   & 0.0 & 0.0 & 1.0 \\ 
\hline
\multicolumn{9}{c}{Problem RP2 (Configuration 2 in \cite{KurganovTadmor2002})} \\
\hline
     & \multicolumn{4}{c|}{$x \leq 0$} & \multicolumn{4}{|c}{$x>0$} \\
\hline
     & $\rho$ & $u$ & $v$  & p & $\rho$ & $u$ & $v$  & p  \\ 
\hline
$y > 0$    & 0.5197 & -0.7259 & 0.0     & 0.4 & 1.0    & 0.0 &  0.0    & 1.0 \\ 
$y \leq 0$ & 1.0    & -0.7259 & -0.7259 & 1.0 & 0.5197 & 0.0 & -0.7259 & 0.4 \\ 
\hline
\multicolumn{9}{c}{Problem RP3 (Configuration 4 in \cite{KurganovTadmor2002})} \\
\hline
     & \multicolumn{4}{c|}{$x \leq 0$} & \multicolumn{4}{|c}{$x>0$} \\
\hline
     & $\rho$ & $u$ & $v$  & p & $\rho$ & $u$ & $v$  & p  \\ 
\hline
$y > 0$    & 0.5065 &  0.8939 & 0.0     & 0.35 & 1.1    & 0.0 &  0.0    & 1.1  \\ 
$y \leq 0$ & 1.1    &  0.8939 & 0.8939  & 1.1  & 0.5065 & 0.0 &  0.8939 & 0.35 \\ 
\hline 
\end{tabular} 
\end{center}
\label{tab.rp2d.ic}
\end{table} 

\begin{figure}[!htbp]
\begin{center}
\begin{tabular}{cc} 
\includegraphics[width=0.47\textwidth]{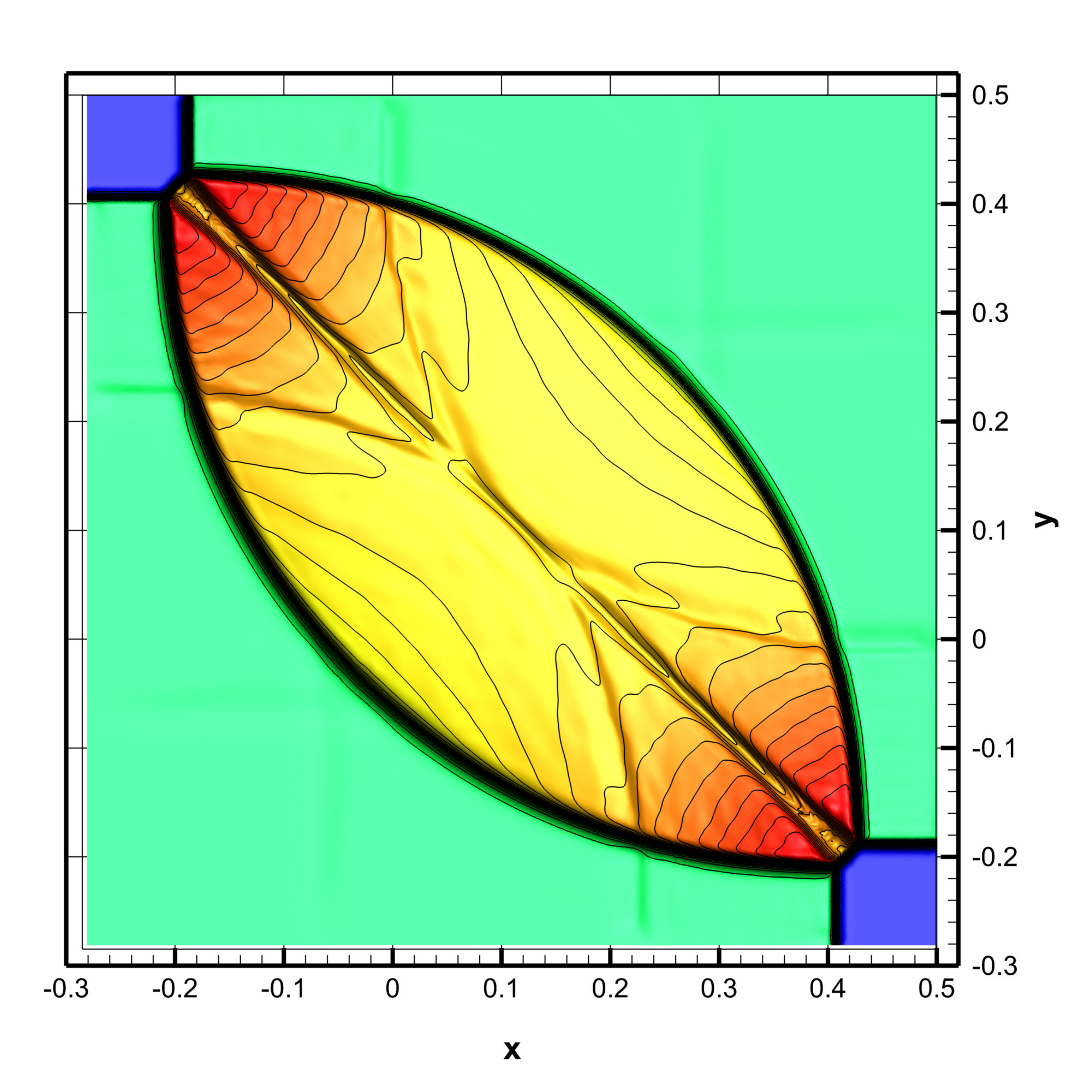}  &           
\includegraphics[width=0.47\textwidth]{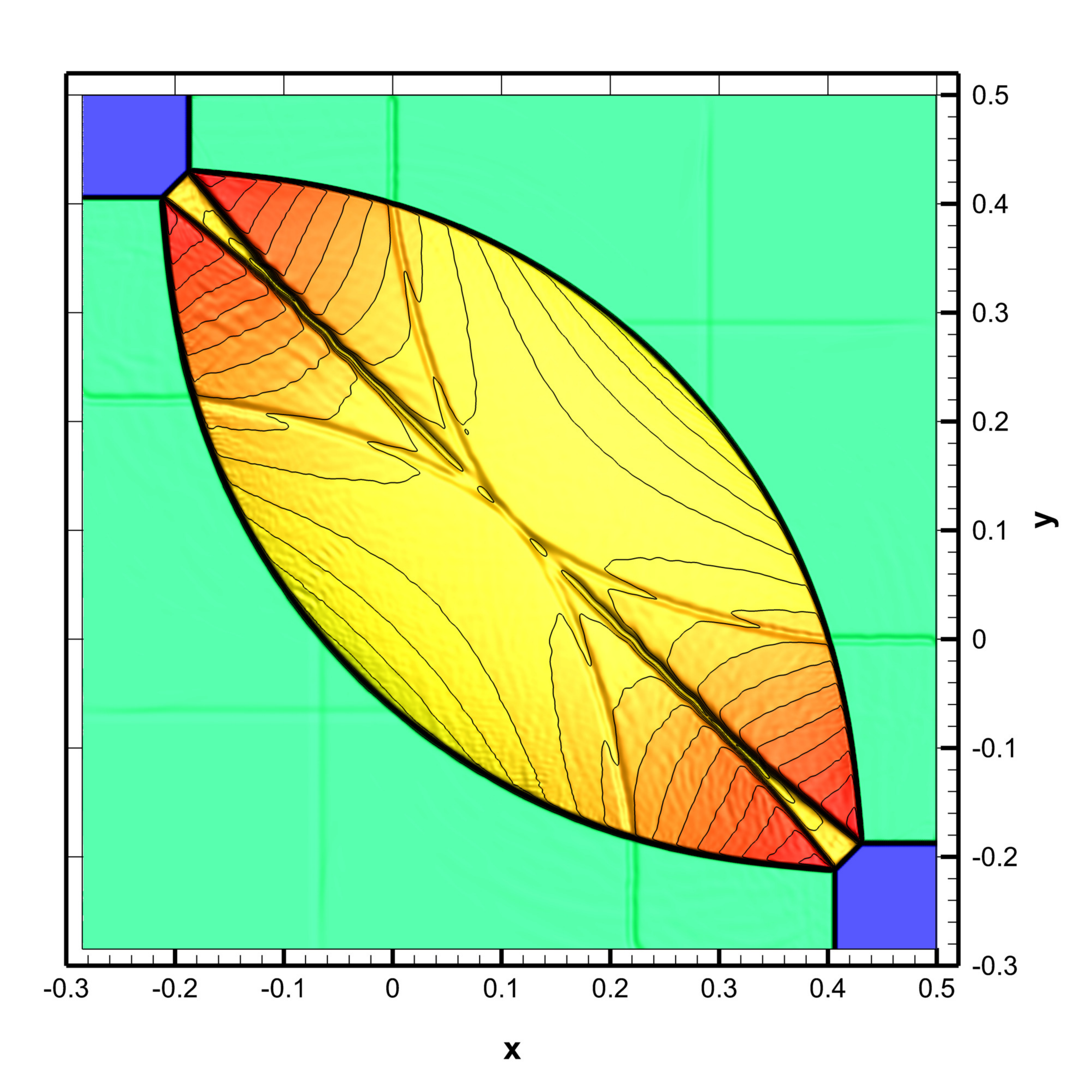} \\   
\multicolumn{2}{c}{
\includegraphics[width=0.70\textwidth]{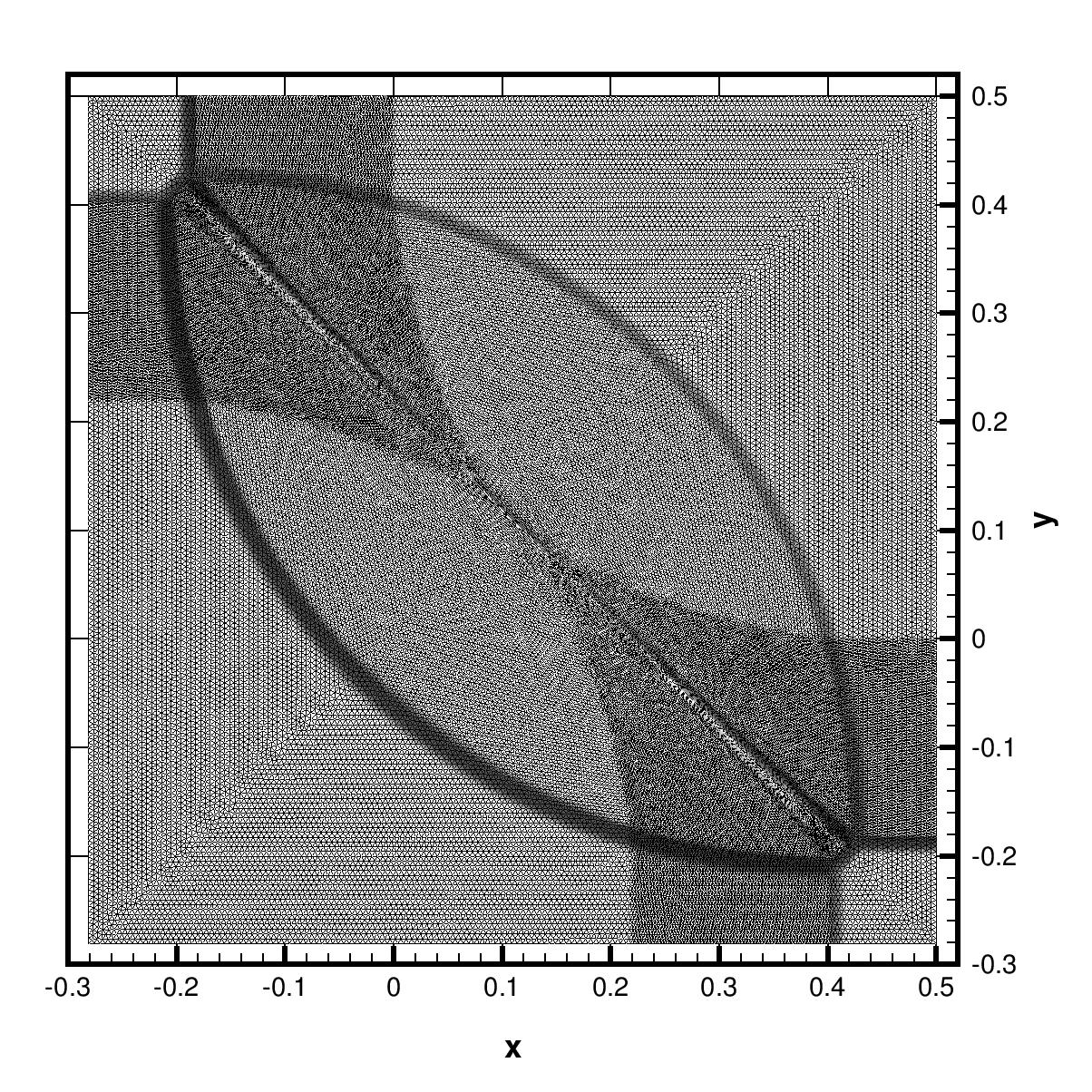}    
}
\end{tabular} 
\caption{Top left: 51 equidistant density contour lines from $\rho=0.53$ to $\rho=1.9$ for the numerical solution obtained with a third order unstructured Lagrangian one--step WENO finite volume scheme 
for the two--dimensional Riemann problem RP3 at $t=0.245$. Top right: Same contour lines for the reference solution. Bottom: computational mesh of the Lagrangian scheme at the final time $t=0.245$.} 
\label{fig.rp3}
\end{center}
\end{figure}

\section{Conclusions} 
\label{sec.conclusions} 
We have developed a new high order two--dimensional Arbitrary--Lagrangian-Eulerian one--step WENO finite volume scheme on unstructured triangular meshes. The algorithm works for general hyperbolic balance laws with non stiff algebraic source terms. Several smooth and non--smooth test problems with shock waves, contact waves, shear waves and rarefactions have been simulated. The results have been compared with exact or numerical reference solutions 
in order to validate our approach. The algorithm was found to work properly in all cases and to be robust and accurate, too. The accuracy has been verified empirically via numerical convergence studies on a smooth test case with exact solution. Further work will concern the improvement of the presented algorithm in order to deal with stiff algebraic source terms, which is straightforward, see the work presented in \cite{Dumbser2012}. Moreover, we plan 
to extend our scheme to three space dimensions in the more general framework of the new $P_{N}P_{M}$ method proposed in \cite{Dumbser20088209}, where one can deal with pure finite volume or pure finite element methods, or with a hybridization of both. Further work will also consist in a generalization to moving \textit{curved} meshes as well as to hyperbolic PDE with \textit{non--conservative} products.   

\section*{Acknowledgments}
The presented research has been financed by the European Research Council (ERC) under the European Union's Seventh Framework 
Programme (FP7/2007-2013) with the research project \textit{STiMulUs}, ERC Grant agreement no. 278267. 

\bibliography{Lagrange2D}
\bibliographystyle{plain}

\end{document}